\documentclass[11pt,leqno]{amsart}
\usepackage{epsfig}
\usepackage{amssymb}
\usepackage{amscd}
\usepackage[matrix,arrow]{xy}

\newtheorem{theo}{Theorem}[section]
\newtheorem{lemma}[theo]{Lemma}
\newtheorem{defi}[theo]{Definition}
\newtheorem{prop}[theo]{Proposition}

\newtheorem{cor}[theo]{Corollary}
\newtheorem{remark}[theo]{Remark}

\newtheorem{example}[theo]{Example}
\numberwithin{equation}{section}

\def\Z{\mathbb{Z}}

\def\coh{\operatorname{coh}}
\def\Qcoh{\operatorname{Qcoh}}

\def\X{{\mathcal{X}}}

\def\bR{{\mathbf R}}
\def\bL{{\mathbf L}}

\def\PP{{\mathbb P}}

\def\pre-tr{\operatorname{pre-tr}}
\def\h{\operatorname{h}}
\def\Hom{\operatorname{Hom}}
\def\End{\operatorname{End}}
\def\gr{\operatorname{gr}}

\newcommand{\cQ}{{\mathcal Q}}
\newcommand{\cF}{{\mathcal F}}
\newcommand{\cG}{{\mathcal G}}
\newcommand{\cO}{{\mathcal O}}
\newcommand{\cP}{{\mathcal P}}
\newcommand{\cL}{{\mathcal L}}
\newcommand{\cM}{{\mathcal M}}
\newcommand{\cN}{{\mathcal N}}

\newcommand{\cA}{{\mathcal A}}
\newcommand{\cB}{{\mathcal B}}

\newcommand{\cC}{{\mathcal C}}
\newcommand{\cE}{{\mathcal E}}

\newcommand{\cR}{{\mathcal R}}

\newcommand{\cH}{{\mathcal H}}

\newcommand{\cl}{\operatorname{cl}}

\newcommand{\DG}{\operatorname{DG}}
\newcommand{\Fun}{\operatorname{Fun}}
\newcommand{\Tot}{\operatorname{Tot}}

\newcommand{\Def}{\operatorname{Def}}
\newcommand{\Perf}{\operatorname{Perf}}

\newcommand{\Ker}{\operatorname{Ker}}
\newcommand{\im}{\operatorname{Im}}

\newcommand{\Ext}{\operatorname{Ext}}
\newcommand{\Tor}{\operatorname{Tor}}

\newcommand{\Aut}{\operatorname{Aut}}

\newcommand{\Spec}{\operatorname{Spec}\,}

\newcommand{\dgart}{\operatorname{dgart}}
\newcommand{\art}{\operatorname{art}}

\newcommand{\coDef}{\operatorname{coDef}}
\newcommand{\cart}{\operatorname{cart}}
\newcommand{\Alg}{\operatorname{Alg}}

\newcommand{\Ho}{\operatorname{Ho}}

\newcommand{\id}{\operatorname{id}}

\newcommand{\DEF}{\operatorname{DEF}}
\newcommand{\coDEF}{\operatorname{coDEF}}

\newcommand{\dg}{\operatorname{dg}}
\newcommand{\coker}{\operatorname{coker}}
\newcommand{\Mod}{\operatorname{Mod}}
\newcommand{\Set}{\operatorname{Sets}}
\newcommand{\Tors}{\operatorname{Tors}}

\newcommand{\QGr}{\operatorname{QGr}}
\newcommand{\QMod}{\operatorname{QMod}}
\newcommand{\qmod}{\operatorname{qmod}}
\newcommand{\Proj}{\operatorname{Proj}}
\newcommand{\QCoh}{\operatorname{QCoh}}
\newcommand{\Coh}{\operatorname{Coh}}
\newcommand{\NGr}{\operatorname{NGr}}
\newcommand{\Gr}{\operatorname{Gr}}
\newcommand{\rank}{\operatorname{rank}}
\newcommand{\pr}{\operatorname{pr}}

\usepackage{epsf}
\usepackage{amscd}

\newcommand{\T}{\mathcal{T}}

\newcommand{\m}{\mathfrak{m}}

\newcommand{\one}{\mathbf{1}}

\title[Deformation theory of objects  in homotopy and derived categories III]{Deformation theory of objects  in homotopy and derived categories III: abelian categories}

\author{Alexander I. Efimov}
\address{Department of Mechanics and Mathematics, Moscow State University, Moscow,
Russia} \email{efimov@mccme.ru}
\author{Valery A.~Lunts}
\address{Department of Mathematics, Indiana University,
Bloomington, IN 47405, USA} \email{vlunts@indiana.edu}
\author{Dmitri O.~Orlov}
\address{Algebra Section, Steklov Mathematical Institute, 8 Gubkina St.,  Moscow, 119991 Russia}
\email{orlov@mi.ras.ru}

\thanks{The first named author was partially supported by the  NSh grant 1983.2008.1, and by
the Moebius Contest Foundation for Young Scientists. The second
named author was partially supported by the NSA grant
H98230-05-1-0050 and CRDF grant RUM1-2661-MO-05.  The third named
author was partially supported by  RFFI grant 05-01-01034, INTAS
grant 05-1000008-8118, and  NSh grant 9969.2006.1}

\begin{document}

\begin{abstract} This is the third paper in a series. In part I we
developed a deformation theory of objects in homotopy and derived
categories of DG categories. Here we show how this theory can be
used to study deformations of objects in homotopy and derived
categories of abelian categories. Then we consider examples from
(noncommutative) algebraic geometry. In particular, we study
noncommutative Grassmanians that are  true noncommutative moduli
spaces of structure sheaves of projective subspaces in  projective
spaces.
\end{abstract}

\maketitle

\tableofcontents

\section{Introduction}

In our paper \cite{ELOI} we developed the deformation theory of a
right DG module over a DG category $\cA $ in the corresponding
homotopy and derived categories. In the subsequent paper
\cite{ELOII} we proved pro-representability of the corresponding
deformation pseudo-functors. In this paper we would like to show
how to apply this theory to deformations of complexes over {\it
abelian} categories (in the corresponding homotopy and derived
categories).

In the second part of the paper we discuss the example of
complexes of (quasi-) coherent sheaves on a scheme. Then we give
examples when when our pro-representabity theorems in \cite{ELOII}
can be applied to this geometric situation.

The third part is devoted to the example of a global noncommutative
moduli space of objects in derived categories: noncommutative
Grassmanians $\NGr(m,V).$ The noncommutative scheme $\NGr(m,V)$ is a
true noncommutative moduli space of structure sheaves
$\cO_{\PP(W)}\in D^b_{coh}(\PP(V)),$ where $W\subset V$ are vector
subspaces of dimension $\dim W=m.$ Namely, it satisfies the
following properties:

1) There is a natural fully faithful functor $\Phi$ from the
category of perfect complexes
$\Perf(\NGr(m,V))$ (Definition
\ref{perfect}) to $D^b_{coh}(\PP(V)).$ Its image is the double
orthogonal to the family of objects $\cO_{\PP(W)},$ i.e. the full
subcategory generated by objects
$\cO_{\PP(V)}(m-n),\dots,\cO_{\PP(V)}(-1),\cO_{\PP(V)}.$ This is
Corollary \ref{T^m,V} below;

2) There is a $k$-point $x_W\in \NGr(m,V)(k)=X_{\cA^{m,V}}(k)$
(see Section \ref{k-points} below) for each subspace $W\subset V$ of
dimension $\dim W=m.$ Further, $(x_W)_*(k)$ lies in
$\Perf(\NGr(m,V))$ and $\Phi(x_*(k))\cong \cO_{\PP(W)}.$ This is a
part of Theorem \ref{all_k-points} below.

3) The completion of the local ring of the $k$-point $x_W$ (see
Section \ref{defi_local}) is isomorphic to $H^0(\hat{S})^{op},$
where $\hat{S}$ is dual to the bar construction of the minimal
$A_{\infty}$-structure on $\Ext^{\cdot}(\cO_{\PP(W)},\cO_{\PP(W)})$
(Theorem \ref{local_rings}).  It can be shown that the DG algebra
$\bR\Hom^{\cdot}(\cO_{\PP(W)},\cO_{\PP(W)})$ is formal and the graded
algebra $\Ext^{\cdot}(\cO_{\PP(W)},\cO_{\PP(W)})$ is quadratic
Koszul, and hence the projection $\hat{S}\to H^0(\hat{S})$ is a
quasi-isomorphism. Hence, the moduli space is not a DG space but
just noncommutative space.

We do not have a moduli functor of our family of objects
$\cO_{\PP(W)},$ which should be defined on the category of
noncommutative affine schemes. However, the properties 1), 2) and 3)
suggest that $\NGr(m,V)$ is a true moduli space of this family of
objects, in our context of deformations of objects in derived
categories.

It is remarkable that there is a natural morphism from the
commutative Grassmanian $\Gr(m,V)$ to noncommutative one
$\NGr(m,V).$ Moreover, the functor $\Phi:\Perf(\NGr(m,V))\to
D^b_{coh}(\PP(V))$ above coincides with $\bL f_{1,m,V}^*,$ where
$f_{1,m,V}:\PP(V)\to \NGr(m,V)$ is a natural morphism. Both these
statements are parts of Proposition \ref{Gr_to_NGr} below.

Section \ref{Z-algebras} contains some preliminaries on
$\Z$-algebras and the associated noncommutative schemes (or
stacks) $\Proj(\cA)$ regarded as an abelian category of
quasi-coherent sheaves together with a structure sheaf.

In section \ref{Grassmans} we define the noncommutative
Grassmanians as $\Proj$ of certain $\Z\text{-}$algebras.

In section \ref{derived_on_NGr} we describe the derived categories
of noncommutative Grassmanians (Theorem \ref{D(NGr)}). This is an
application of the more general result for geometric
$\Z\text{-}$algebras (Theorem \ref{D(Geom)}) which originally
appeared in \cite{BP}.

In  section \ref{k-points} we make an attempt to relate two
different approaches to noncommutative geometry. Namely, we
associate to each (positively oriented) $\Z\text{-}$algebra a
presheaf of groupoids $X_{\cA}$ on the category $\Alg_k^{op}$ dual
to the category of associative unital $k$-algebras. The groupoid
$X_{\cA}(B)$ should be thought of as a groupoid of maps $Sp(B)\to
\Proj(\cA).$ We compare our definition with maps between commutative
schemes (Proposition \ref{comm_situation}). Then we describe the
k-points of noncommutative Grassmanians (Theorem
\ref{all_k-points}).

In the last section \ref{local_structure}, for any presheaf $X$ of
sets on $\Alg_k^{op},$ and its $k\text{-}$point $x\in X(k),$ we
define the notion of a completion of local ring $\widehat{\cO_x},$
which can not exist a priori. Then we prove that in the case of the
noncommutative Grassmanian $\NGr(m,V)$ and the $k\text{-}$point
$x_W$ corresponding to the subspace $W\subset V$ of dimension $m,$
the completion $\widehat{\cO_x}$ exists and is isomorphic to
$H^0(\hat{S})^{op} (\cong H^0(\hat{S}))$ in the above notation.

We freely use the notation and results of \cite{ELOI} and
\cite{ELOII}. The reference to \cite{ELOI} or \cite{ELOII} appears
in the form I, Theorem ... , or II, Theorem ... respectively.

\part{Deformations of objects in homotopy and derived categories of
abelian categories}

Let $\cM$ be small a $k$-linear abelian category. Denote by
$C(\cM),$ $\cH(\cM),$ $D(\cM)$ the category of complexes over $\cM,$
its homotopy category and its derived category respectively. We will
also consider the usual categories $C^b(\cM),$ $C^{\pm}(\cM)$ of
bounded (resp. bounded above, below) complexes and the categories
$\cH ^b(\cM),$ $D^b(\cM),$ $\cH^{\pm}(\cM),$ $D^{\pm}(\cM)$ of
cohomologically bounded (resp. bounded below, above) complexes.
Given $E\in C(\cM)$ and an artinian DG algebra $\cR$ there are
natural notions of homotopy and derived $\cR$-deformations (and
$\cR$-co-deformations) of $E.$ We start by defining this deformation
theory and then show (under some assumptions) how it can be
interpreted as a deformation theory of a DG module over an
appropriate DG category $\cA.$  This interpretation allows us to
translate the previous results obtained in the DG context to the
case of $C(\cM).$ Our point of view is that the deformation theory
developed in \cite{ELOI}, \cite{ELOII} in the language of DG modules
is more flexible. Hence for example in the context of abelian
categories we omit the notion of pseudo-functors $\DEF $ and
$\coDEF$ from the $2\text{-}\dgart$ to ${\bf Gpd}.$

\section{Categories $C_{\cR}(\cM),$ $C_{\cR}^{dg}(\cM),$ $\cH _{\cR}(\cM),$
$D_{\cR}(\cM)$}

The category of complexes over $\cM$ is also naturally a DG category
with the $\Hom$-complexes being the usual complexes of morphisms
between objects in $C(\cM).$ We denote this DG category by
$C^{dg}(\cM).$ Then $\cH(\cM)$ is simply $\Ho(C^{dg}(\cM))$ and the
category $D(\cM)$ is obtained from $\Ho(C^{dg}(\cM))$ by inverting
quasi-isomorphisms. Notice that this is NOT the same as
$D(C^{dg}(\cM))$ as defined in I, Section 3.1.

\begin{defi}\label{D_R_etc} Let $\cR$ be an artinian DG algebra. A right
$\cR$-complex over $\cM$ (or simply an $\cR ^{op}$-complex) is an
object $S\in C(\cM)$ together with a homomorpism of DG algebras $\cR
^{op}\to \Hom^{\cdot} (S,S).$ This is the same as a DG functor from
the DG category $\cR ^{op}$ (with one object) to the DG category
$C^{dg}(\cM).$ Thus $\cR ^{op}$-complexes over $\cM$ naturally form
a DG category $\Fun_{\DG}(\cR ^{op},C^{dg}(\cM))$ which we denote by
$C^{dg}_{\cR}(\cM).$ If in the category $C^{dg}_{\cR}(\cM)$ we only
consider morphisms which are degree zero cycles (i.e. DG
transformations between DG functors), then we obtain an abelian
category, which we denote by $C_{\cR}(\cM).$ The homotopy category
$\Ho(C^{dg}_{\cR}(\cM))$ is denoted by $\cH _{\cR}(\cM).$ If we
invert quasi-isomorphisms in $\cH _{\cR}(\cM)$ we obtain the derived
category $D_{\cR}(\cM).$
\end{defi}

The categories $\cH _{\cR}(\cM)$and $D_{\cR}(\cM)$ are naturally
triangulated. We will also consider the obvious full DG
subcategories $C^b_{\cR}(\cM), C^{\pm}_{\cR}(\cM)\subset
C_{\cR}^{dg}(\cM)$ and the full triangulated subcategories $\cH^b
_{\cR}(\cM), \cH ^{\pm}_{\cR}(\cM)\subset \cH _{\cR}(\cM)$;
$D^b_{\cR}(\cM), D^{\pm}_{\cR}(\cM)\subset D_{\cR}(\cM).$

\begin{remark}\label{abelian_vs_dg} \rm{Consider $\cM$ as a DG category (with all morphisms
being of degree zero) and let $\cR$ be an artinian DG algebra.
Notice that an $\cR ^{op}$-complex $S\in C(\cM)$ defines (by
Yoneda) a DG-module over the DG category $\cM _{\cR}^{op}=\cM
^{op}\otimes \cR ^{op}$ (3.1, 3.3 in Part I), i.e. there is a full
and faithful embedding of DG categories
$$h^\bullet _{\cR}: C^{dg}_{\cR}(\cM)\hookrightarrow \cM ^{op}_{\cR}\text{-mod}.$$
Using this embedding we could directly apply our machinery in Part I
to obtain a deformation theory of objects in $C(\cM).$ This
deformation theory however would not always give the right answer
(in case of derived deformations). Our point is that there exists a
natural independent deformation theory for complexes over abelian
categories which we define in the next section. Eventually we will
compare this theory to deformations of DG-modules as in Part I.}
\end{remark}

The next lemma is a repetition of I, Lemma 3.19 in our context.

\begin{lemma}\label{truncation} Assume that $\cR \in \dgart _-.$ Then there exist
truncation functors in $D_{\cR}(\cM)$: for every
$\cR^{op}$-complex $S$ there exists a short exact sequence of $\cR
^{op}$-complexes
$$\tau _{<0}S \to S \to \tau _{\geq 0}S,$$
where $H^i(\tau _{<0}S)=0$ if $i\geq 0$ and $H^i(\tau _{\geq 0}S)=0$
for $i<0.$
\end{lemma}

\begin{proof} Indeed, put
$$\tau _{<0}S:=\oplus_{i<0}S^i\oplus d(S^{-1}).$$
\end{proof}

The next definition is the analogue of I, Definition 3.8.

\begin{defi}\label{R-free} Let $\cR \in \dgart.$ An $\cR ^{op}$-complex $S$ is called
graded $\cR$-free (resp. graded $\cR$-cofree) if there exists
$M\in C^{dg}(\cM)$ and an isomorphism of graded objects
(forgetting the differential) in $C^{dg}_{\cR}(\cM)$ $M\otimes
\cR\simeq S$ (resp. $M\otimes \cR ^*\simeq S$).
\end{defi}

\begin{prop}\label{phi^*_*^!} A homomorphism $\phi :\cR \to \cQ$ of artinian DG algebras induces
DG functors
$$\phi ^*:C^{\dg}_{\cR}(\cM)\to C^{\dg}_{\cQ}(\cM),\quad \phi
_*:C^{\dg}_{\cQ}(\cM) \to C^{\dg}_{\cR}(\cM),\quad \phi
^!:C^{\dg}_{\cR}(\cM) \to C^{\dg}_{\cQ}(\cM).$$ The DG functors
$(\phi ^*,\phi _*)$ and $(\phi _*,\phi^!)$ are adjoint. That is for
$S\in C^{\dg}_{\cR}(\cM)$ and $T\in C^{\dg}_{\cQ}(\cM)$ there are
functorial isomorphisms of complexes
$$\Hom^{\cdot} (\phi ^*S,T)=\Hom^{\cdot} (S,\phi _*T),\quad \Hom^{\cdot} (\phi _*T,S)=\Hom^{\cdot}
(T,\phi ^!S).$$
  We denote by the
same symbols the induced functors between the abelian categories
$C_{\cR}(\cM),$ $C_{\cQ}(\cM)$ and the homotopy categories
$\cH_{\cR}(\cM),$ $\cH_{\cQ}(\cM).$   These induced functors are
also adjoint.
\end{prop}

\begin{proof} The categories $\cM$ and $C(\cM)$ are abelian and as such have all
finite limits and colimits. Let $S\in C^{dg}_{\cR}(\cM).$ We put
$$\phi ^*(S)=S\otimes _{\cR}Q,\quad  \phi ^!(S)=\Hom^{\cdot} _{\cR ^{op}}(Q,S).$$
That is $\phi ^*(S)$ is defined as a colimit of a finite (since
$\dim \cR <\infty$) diagram involving the object $S\otimes _k\cQ$;
and $\phi ^!(S)$ is a limit of a finite (since $\dim \cR<\infty$)
diagram involving the object $S\otimes _kQ^*$ (since $\dim
Q<\infty$). So these objects are well defined. The DG functor $\phi
_*$ is simply the restriction of scalars.

For each $M\in C^{dg}_{\cQ}(\cM),$ $N\in C^{dg}_{\cR}(\cM)$ we
have natural functorial closed morphisms of degree zero
$$\eta_1(M):\phi^*\phi_*(M)=\phi_*(M)\otimes_{\cR}\cQ\to M,\quad
\eta_2(N):N\to \phi_*(N\otimes_{\cR}\cQ)=\phi_*\phi^*(N),$$
induced by the structure morphism $M\otimes\cQ\to M,$ and by the
inclusion $N\to N\otimes \cQ.$ They give rise to the morphisms of
functors
$$\eta_1:\phi^*\phi_*\to id,\quad \eta_2:id\to \phi_*\phi^*.$$ It is also clear that the
compositions
$$\phi^*\stackrel{\phi^*(\eta_2)}{\to} \phi^*\phi_*\phi^*\stackrel{\eta_1(\phi^*)}\to \phi^*,$$
$$\phi_*\stackrel{\eta_2(\phi_*)}{\to} \phi_*\phi^*\phi_*\stackrel{\phi_*(\eta_1)}\to \phi_*$$
are equal to identity morphisms. Hence, the pair of DG functors
$(\phi^*,\phi_*)$ is adjoint. The same arguments hold for abelian
and homotopy categories.

The adjunction $(\phi_*,\phi^!)$ is proved in the same way.
\end{proof}

\begin{example}\label{i_and_p} \rm{Let $i:\cR \to k$ be the augmentation map and
$p:k\to \cR$ be the obvious inclusion. Then we obtain the
corresponding functors $i^*,i_*,i^!,p^*,p_*,p^!.$}
\end{example}

\begin{defi}\label{h-proj_inj} An object $S\in C^{\dg}_{\cR}(\cM)$ is called
h-projective (resp. h-injective) if for every acyclic $T\in
C^{\dg}_{\cR}(\cM)$ the complex $\Hom^{\cdot} (S,T)$ (resp.
$\Hom^{\cdot} (T,S)$) is acyclic.
\end{defi}

\begin{remark}\label{proj_dir_sums} \rm{Note that the collection of h-projectives (resp.
h-injectives) is closed under arbitrary (existing in
$C^{\dg}_{\cR}(\cM)$) direct sums (resp. direct products).}
\end{remark}

\begin{cor}\label{phi^*(proj)} The DG functor $\phi ^*$ (resp. $\phi ^!$) preserves
h-projectives (resp. h-injectives).
\end{cor}

\begin{proof} This follows from the adjunctions $(\phi ^*,\phi _*),$
$(\phi _*,\phi ^!)$ and the fact that $\phi _*$ preserves acyclic
complexes.
\end{proof}

\begin{prop}\label{enough} Fix $\cR \in \dgart _-.$

a) Assume that $\cM$ has enough projectives. Then for every $S\in
C_{\cR}^-(\cM)$ there exists an h-projective $P\in C_{\cR}^-(\cM)$
and a quasi-isomorphism $P\to S.$ We may choose $P$ to be graded
$\cR$-free.

b)  Assume that $\cM$ has enough injectives. Then for every $T\in
C_{\cR}^+(\cM)$ there exists an h-injective $I\in C_{\cR}^+(\cM)$
and a quasi-isomorphism $T\to I.$ We may choose $I$ to be graded
$\cR$-cofree.
\end{prop}

\begin{proof} a) Fix $S\in
C_{\cR}^-(\cM)$ and assume  that $S^i=0$ for $i>i_0.$ Since $\cM$
has enough projectives we can find (by a standard construction) an
h-projective $Q_0\in C^-(\cM)$ and a surjective quasi-isomorphism
$\epsilon ^\prime:Q_0\to p_*S.$ We may and will assume that each
$Q_0^j\in \cM$ is projective. Moreover we may and will assume that
$Q_0^i=0$ for $i>i_0.$ By adjunction we obtain a surjective morphism
in $C^-_{\cR}(\cM)$
$$\epsilon :P_0=p^*Q_0\to S,$$
which is also surjective on cohomology. Denote $K=\Ker (\epsilon).$
Note that since $\cR \in \dgart_-$ we have $K^i=0$ for $i>i_0.$ Now
repeat the procedure with $K$ instead of $S.$ Finally we obtain an
exact sequence $$\cdots\to P_{-2}\to P_{-1}\to
P_{0}\stackrel{\epsilon}{\to}S\to 0$$ in $C_{\cR}^{-}(\cM)$ with the
following properties:

i) $P_{-n}$ is h-projective for every $n.$

ii) $P_{-n}^i=0$ for $i>i_0$ and all $n.$

iii) The complex $$\cdots\to H^{\cdot}(P_{-1})\to H^{\cdot}(P_0)\to
H^{\cdot}(S)\to 0$$ is exact.

Denote by $P$ the total complex
$$P=\Tot (\cdots\to P_{-1}\to P_0).$$
Note that as a graded object
$$P=\oplus _{n\geq 0}P_{-n}[n],$$
so that in each degree $i$ the contribution to $P^i$ comes from
finitely many $P_{-n}$'s. Thus $P$ is a well defined object in
$C_{\cR}^-(\cM)$ (we do not assume that infinite direct sums exist
in $\cM$). Because of the property iii) above the morphism
$$\epsilon :P\to S$$
is a quasi-isomorphism. It remains to show that $P$ is h-projective.

We have the standard increasing filtration by
$\cR^{op}\text{-}$subcomplexes $$F_nP=\Tot (P_{-n}\to \dots \to
P_{-1}\to P_0).$$ This filtration satisfies the following
properties:

(a) $P=\bigcup\limits_{n\geq 0} F_n P;$

(b) each quotient $F_n P/F_{n-1}P$ is h-projective;

(c) each inclusion of graded $\cR^{op}\text{-}$objects
$(F_{n-1}P)^{gr}\hookrightarrow (F_n P)^{gr}$ splits.

It follows that $P$ is h-projective.

b) The proof is very similar to that of a), but we present it anyway
for completeness. Fix $T\in C_{\cR}^+(\cM),$ say $T^i=0$ for
$i<i_0.$ Since $\cM$ has enough injectives we can find an
h-injective $J_0\in C^+(\cM)$ and an injective quasi-isomorphism
$s^\prime :p_*T\to J_0.$ We may and will assume that $J_0$ consists
of objects which are injective in $\cM.$ Moreover we may and will
assume that $J_0^i=0$ for $i<i_0.$ By adjunction we obtain an
injective morphism of objects in $C^+_{\cR}(\cM)$
$$s:T\to I_0:=p^!J_0,$$
which is also injective on cohomology. Denote $L=\coker(s).$ Note
that since $\cR \in \dgart_-$ we have $L^i=0$ for $i<i_0.$ Now
repeat the procedure with $L$ instead of $T.$ Finally we obtain an
exact sequence
$$0\to T \stackrel{s}{\to} I_0 \to I_1\to ...$$
in $C_{\cR}^+(\cM)$ with the following properties:

i') $I_n$ is h-injective for all $n.$

ii') $I_n^i=0$ for $i<i_0$ and all $n.$

iii') The complex
$$0\to H^{\cdot}(T)\to H^{\cdot}(I_0)\to H^{\cdot}(I_1)\to ...$$
is exact.

Denote by $I$ the total complex
$$I:=\Tot(I_0\to I_1\to ...).$$
Note that as a graded object
$$I=\oplus _{n\geq 0}I_n[-n],$$
so that in each degree $i$ the contribution to $I^i$ comes from
finitely many $I_n$'s. Thus $I$ is a well defined object in
$C^+_{\cR}(\cM).$ Because of property iii') above the morphism
$s:T\to I$ is a quasi-isomorphism. It remains to show that $I$ is
h-injective.

We have the standard decreasing filtration by
$\cR^{op}\text{-}$subcomplexes $$F_nI=\Tot (I_n\to I_{n+1}\to
\dots).$$ This filtration satisfies the following properties:

(a) $I=\lim\limits_{\leftarrow} I/F_n I;$

(b) each quotient $F_n P/F_{n+1}P$ is h-injective;

(c) each inclusion of graded $\cR^{op}\text{-}$objects
$(F_{n+1}I)^{gr}\hookrightarrow (F_n P)^{gr}$ splits.

It follows that $I$ is h-injective.

Proposition is proved.
\end{proof}

Using the last proposition we can define derived functors of the
functors $\phi ^*$ and $\phi ^!.$ Namely assume that $\cM$ has
enough projectives (resp. injectives). Then given a homomorphism
$\phi :\cR \to \cQ$ of artinian (non-positive) DG algebras we define
the functor $\bL \phi ^*:D^-_{\cR}(\cM)\to D^-_{\cQ}(\cM)$ (resp.
 $\bR \phi ^!:D^+_{\cR}(\cM)\to D^+_{\cQ}(\cM)$) using h-projectives
 (resp. h-injectives) in the usual way. Notice that the functor
 $\phi _*$ is exact, hence it extends trivially to $\phi _*:D^{\pm}_{\cQ}(\cM)\to
 D^{\pm}_{\cR}(\cM).$ The functors $(\bL \phi ^*,\phi _*)$ and
 $(\phi _*,\bR \phi ^!)$ are adjoint.

\begin{remark}\label{RHom} \rm{Let $\cR \in \dgart _-$ and $M,N \in
C^{dg}_{\cR}(\cM).$ Assume that i) $\cM$ has enough projectives and
$M\in D^-_{\cR}(\cM)$ or ii) $\cM$ has enough injectives and $N\in
D^+_{\cR}(\cM).$ Then we can define the complex $\bR \Hom^{\cdot}
(M,N)$ and hence the vector spaces $\Ext ^i(M,N).$ Namely, by
Proposition \ref{enough} in the first case we may replace $M$ by a
quasi-isomorphic h-projective $P$ and in the second case we may
replace $N$ by a quasi-isomorphic h-injective $I.$ Then define $\bR
\Hom (M,N)$ as either $ \Hom^{\cdot} (P,N)$ or $\Hom^{\cdot}
(M,I).$}
\end{remark}

\begin{prop}\label{invariance} Let $\phi :\cR \to \cQ$ be a morphism of artinian DG
algebras which is a quasi-isomorphism.

a) Assume that $\cM$ has enough projectives. Then the functor
$\bL\phi ^*:D^-_{\cR}(\cM)\to D^-_{\cQ}(\cM)$ is an equivalence of
categories;

b) Assume that $\cM$ has enough injectives. Then the functor
$\bR\phi ^!:D^+_{\cR}(\cM)\to D^+_{\cQ}(\cM)$ is an equivalence of
categories.
\end{prop}

\begin{proof} a) It suffices to prove that for each h-projective
$P\in C^{dg}_{\cR}(\cM),$ $P'\in C^{dg}_{\cQ}(\cM),$ the
adjunction morphisms $$\eta_1(P'):\phi^*\phi_*(P')\to P',\quad
\eta_2(P):P\to \phi_*\phi^*(P)$$ are quasi-isomorphisms. As
complexes in $C(\cM),$ the cones of both morphisms are of the form
$P''\otimes_{\cR}Cone(\cR\to \cQ),$ where $P''$ is $P$ or $P'$
respectively. So, it remains to prove the following

\begin{lemma} If $N$ is an acyclic finite-dimensional
$\cR\text{-}$module and $P\in C^{dg}_{\cR}(\cM)$ is h-projective,
then the complex $P\otimes_{\cR}N$ is null-homotopic.
\end{lemma}
\begin{proof} Let $M$ be an object in $C(\cM).$ Then we have
$$\Hom^{\cdot}(P\otimes_{\cR}N,M)=\Hom^{\cdot}_{\cR^{op}}(P,\Hom^{\cdot}_k(N,M)),$$ and
this complex is acyclic since the $\cR^{op}\text{-}$complex
$\Hom^{\cdot}_k(N,M)$ is acyclic and $P$ is h-projective. Hence,
the complex $P\otimes_{\cR}N$ is null-homotopic.
\end{proof}

b) The proof is similar. It suffices to prove that for each
h-injective $I\in C^{dg}_{\cR}(\cM),$ $I'\in C^{dg}_{\cQ}(\cM),$
the adjunction morphisms $$\eta_1(I):\phi_*\phi^!(I)\to I,\quad
\eta_2(I'):I'\to \phi^!\phi_*(I')$$ are quasi-isomorphisms. As
complexes in $C(\cM),$ the cones of both morphisms are of the form
$\Hom^{\cdot}_{\cR^{op}}(Cone(\cR\to \cQ),I'')$ where $I''$ is $I$
or $I'$ respectively. So, it remains to prove the following

\begin{lemma} If $N$ is an acyclic finite-dimensional
$\cR^{op}\text{-}$module and $I\in C^{dg}_{\cR}(\cM)$ is
h-projective, then the complex $\Hom^{\cdot}_{\cR^{op}}(N,I)$ is
null-homotopic.
\end{lemma}
\begin{proof} Let $M$ be an object in $C(\cM).$ Then we have
$$\Hom^{\cdot}(M, \Hom^{\cdot}_{\cR^{op}}(N,I))=\Hom^{\cdot}_{\cR^{op}}(M,\Hom^{\cdot}_k(N,I)),$$ and
this complex is acyclic since the $\cR^{op}\text{-}$complex
$\Hom^{\cdot}_k(N,M)$ is acyclic, h-injective and hence
null-homotopic. Hence, the complex $\Hom^{\cdot}_{\cR}(N,I)$ is
null-homotopic.
\end{proof}
Proposition is proved.
\end{proof}

 \section{Deformation pseudo-functors}

Let $E$ be an object in $C(\cM).$ As in the DG setting we first
define the homotopy deformation and co-deformation pseudo-functors
$$\Def ^{\h}(E),\  \coDef ^{\h}(E):\dgart \to {\bf
Gpd}.$$ The definitions are copies (in our context) of the
corresponding definitions in Part I.

\begin{defi}\label{def^h} Fix $E \in C(\cM)$ and let $\cR \in \dgart.$ An object
in the groupoid $\Def ^{\h}_{\cR}(E)$ is a pair $(S, \sigma),$ where
$S\in C_{\cR}(\cM)$ and $\sigma :i^*S\to E$ is an isomorphism of
objects in $C(\cM)$ such that the following holds: there exists an
isomorphism of graded objects $\eta :(E\otimes \cR)^{\gr}\to
S^{\gr}$ so that the composition
$$E=i^*(E\otimes \cR)\stackrel{i^*\eta}{\to} i^*S
\stackrel{\sigma}{\to}E$$ is the identity.

Given objects $(S,\sigma ),(S^\prime ,\sigma ^\prime)\in \Def
^{\h}_{\cR}(E),$ a map $f:(S,\sigma )\to (S^\prime ,\sigma ^\prime)$
is an isomorphism $f:S\to S^\prime$ such that $\sigma ^\prime \cdot
i^*f=\sigma.$ An allowable homotopy between maps $f,g$ is a homotopy
$h:f\to g$ such that $i^*h=0.$ We define morphisms in $\Def
^{\h}_{\cR}(E)$ to be classes of maps modulo allowable homotopies.

Note that a homomorphism of artinian DG algebras $\phi :\cR \to \cQ$
induces the functor $\phi ^*:\Def _{\cR}^{\h}(E)\to \Def
_{\cQ}^{\h}(E).$ This defines the pseudo-functor
$$\Def ^{\h}(E):\dgart \to {\bf Gpd}.$$
\end{defi}

We refer to objects of $\Def _{\cR}^{\h} (E)$ as homotopy
$\cR$-deformations of $E.$

\begin{example} \rm{We call $(p^*E,\id)\in \Def _{\cR}^{\h}(E)$ the
trivial $\cR$-deformation of $E.$}
\end{example}

\begin{defi} Denote by $\Def _+^{\h}(E),$ $\Def _-^{\h}(E),$ $\Def
_0^{\h}(E),$ $\Def ^{\h}_{\cl}(E)$ the restrictions of the
pseudo-functor $\Def ^{\h}(E)$ to subcategories $\dgart _+,$ $\dgart
_-,$ $\art,$ $\cart$ respectively.
\end{defi}

Let us give an alternative description of the same deformation
problem. We will define the homotopy {\it co-deformation}
pseudo-functor $\coDef^{\h}(E)$ and (eventually) show that it is
equivalent to $\Def ^{\h}(E).$ The point is that in practice one
should use $\Def ^{\h}(E)$ for a h-projective $E$ and $\coDef
^{\h}(E)$ for a h-injective $E.$

\begin{defi}\label{codef^h} Fix $E\in C(\cM).$ Let $\cR$ be an artinian DG algebra. An object in the groupoid
$\coDef^{\h}_{\cR}(E)$ is a pair $(T, \tau),$ where $T\in
C_{\cR}(\cM)$ and $\tau :E\to i^!T$ is an isomorphism of objects in
$C(\cM)$ so that the following holds: there exists an isomorphism of
graded objects $\delta :T^{\gr}\to (E\otimes \cR ^*)^{\gr}$ such
that the composition
$$E \stackrel{\tau}{\to}i^!T \stackrel{i^!\delta}{\to} i^!(E\otimes \cR ^*)
 =E$$ is the identity.

Given objects $(T,\tau)$ and $(T^\prime,\tau ^\prime)\in \coDef
^{\h}_{\cR}(E)$ a map $g:(T,\tau)\to (T^\prime ,\tau ^\prime)$ is an
isomorphism $f:T\to T^\prime$ such that $i^!f \cdot \tau =\tau
^\prime.$ An allowable homotopy between maps $f,g$ is a homotopy
$h:f\to g$ such that $i^!(h)=0.$ We define morphisms in $\coDef
_{\cR}^{\h}(E)$ to be classes of maps modulo allowable homotopies.

Note that a homomorphism of DG algebras $\phi :\cR \to \cQ$ induces
the functor $\phi ^!:\Def _{\cR}^{\h}(E)\to \Def _{\cQ}^{\h}(E).$
This defines the pseudo-functor
$$\Def ^{\h}(E):\dgart \to {\bf Gpd}.$$
\end{defi}

We refer to objects of $\coDef _{\cR}^{\h} (E)$ as homotopy
$\cR$-co-deformations of $E.$

\begin{example} \rm{For example we can take $T=E\otimes \cR ^*$ with the
differential $d_{E,R^*}:=d_E\otimes 1+1\otimes d_{R^*}$ (and $\tau
=\id$). This we consider as the {\it trivial} $\cR$-co-deformation
of $E.$}
\end{example}

\begin{defi} Denote by $\coDef _+^{\h}(E),$ $\coDef _-^{\h}(E),$ $\coDef
_0^{\h}(E),$ $\coDef _{\cl}^{\h}(E)$ the restrictions of the
pseudo-functor $\coDef ^{\h}(E)$ to subcategories $\dgart _+,$
$\dgart _-,$ $\art,$ $\cart$ respectively.
\end{defi}

\subsection{Derived deformation pseudo-functors}
Likewise we define derived deformation pseudo-functors. In view of
Proposition \ref{enough} above we restrict ourselves to non-positive
artinian DG algebras and consider two cases: enough projectives or
enough injectives in $\cM.$

\begin{defi}\label{def_-} Assume that $\cM$ has enough projectives and fix $E\in D^-(\cM).$
We are going to define a pseudo-functor
$$\Def _-(E):\dgart _-\to {\bf Gpd}$$
of derived deformations of $E.$ Fix an artinian DG algebra $\cR\in
\dgart _-.$ An object of the groupoid $\Def _{\cR}(E)$ is a pair
$(S,\sigma),$ where $S\in D^-_{\cR}(\cM)$ and $\sigma$ is an
isomorphism (in $D^-(\cM)$)
$$\sigma :\bL i^*S\to E.$$
A morphism $f:(S,\sigma)\to (T,\tau)$ between two $\cR$-deformations
of $E$ is an isomorphism (in $D^-_{\cR}(\cM)$) $f:S\to T,$ such that
$$\tau \cdot \bL i^*(f)=\sigma.$$
This defines the groupoid $\Def _{\cR}(E).$ A homomorphism of
artinian DG algebras $\phi:\cR \to \cQ$ induces the functor
$$\bL\phi ^*:\Def _{\cR}(E)\to \Def _{\cQ}(E).$$
Thus we obtain a pseudo-functor
$$\Def _-(E):\dgart _-\to {\bf Gpd}.$$
\end{defi}

We call $\Def _-(E)$ the functor of derived deformations of $E.$

\begin{remark} \rm{A quasi-isomorphism $\phi :\cR\to \cQ$ of artinian
DG algebras induces an equivalence of groupoids
$$\bL\phi ^*:\Def _{\cR}(E)\to \Def _{\cQ}(E).$$ Indeed,
$\bL\phi ^*:D^-_{\cR}(\cM)\to D^-_{\cQ}(\cM)$ is an equivalence of
categories (Lemma \ref{invariance} a)) which commutes with the
functor $\bL i^*.$}
\end{remark}

\begin{remark}\label{quasi-iso_def} \rm{A quasi-isomorphism $\delta:E_1\to E_2$
in $D^-(\cM)$ induces an equivalence of pseudo-functors
$$\delta _*:\Def _-(E_1)\to \Def _-(E_2)$$
by the formula $\delta _*(S,\sigma)=(S,\delta \cdot \sigma).$}
\end{remark}

\begin{defi} Denote by $\Def
_0(E),$ $\Def _{\cl}(E)$ the restrictions of the pseudo-functor
$\Def _-(E)$ to subcategories  $\art,$ $\cart$ respectively.
\end{defi}

Let us now define the pseudo-functor of derived co-deformations.

\begin{defi}\label{codef_-} Assume that $\cM$ has enough injectives and fix $E\in D^+(\cM).$
We are going to define a pseudo-functor
$$\coDef _-(E):\dgart _-\to {\bf Gpd}$$
of derived co-deformations of $E.$ Fix an artinian DG algebra
$\cR\in \dgart _-.$ An object of the groupoid $\coDef _{\cR}(E)$ is
a pair $(S,\sigma),$ where $S\in D^+_{\cR}(\cM)$ and $\sigma$ is an
isomorphism (in $D^+(\cM)$)
$$\sigma :E\to \bR i^!S.$$
A morphism $f:(S,\sigma)\to (T,\tau)$ between two $\cR$-deformations
of $E$ is an isomorphism (in $D^+_{\cR}(\cM)$) $f:S\to T,$ such that
$$ \bR i^!(f)\cdot \sigma=\tau.$$
This defines the groupoid $\coDef _{\cR}(E).$ A homomorphism of
artinian DG algebras $\phi:\cR \to \cQ$ induces the functor
$$\bR\phi ^!:\coDef _{\cR}(E)\to \coDef _{\cQ}(E).$$
Thus we obtain a pseudo-functor
$$\coDef _-(E):\dgart _- \to {\bf Gpd}.$$
\end{defi}

We call $\coDef _-(E)$ the functor of derived co-deformations of
$E.$

\begin{remark} \rm{A quasi-isomorphism $\phi :\cR\to \cQ$ of artinian
DG algebras induces an equivalence of groupoids
$$\bR\phi ^!:\coDef _{\cR}(E)\to \coDef _{\cQ}(E).$$ Indeed,
$\bR\phi ^!:D^+_{\cR}(\cM)\to D^+_{\cQ}(\cM)$ is an equivalence of
categories (Lemma \ref{invariance} b)) which commutes with the
functor $\bR i^!.$}
\end{remark}

\begin{remark}\label{quasi-iso_codef} \rm{A quasi-isomorphism $\delta:E_1\to E_2$
in $D^+(\cM)$ induces an equivalence  of functors
$$\delta ^*:\coDef (E_2)\to \coDef (E_1)$$
by the formula $\delta ^*(S,\sigma)=(S,\sigma \cdot \delta).$}
\end{remark}

\begin{defi} Denote by $\coDef
_0(E),$ $\coDef _{\cl}(E)$ the restrictions of the pseudo-functor
$\coDef _- (E)$ to subcategories  $\art,$ $\cart$ respectively.
\end{defi}

\subsection{Summary of main properties of deformation pseudo-functors}

\begin{prop}\label{def^h=codef^h} Let $E\in C(\cM).$ There exists a natural equivalence of pseudo-functors
$$\Def ^{\h} (E)\to \coDef ^{\h}(E).$$
Consider $E$ as an object in the DG category $C^{dg}(\cM)$ and
denote by $\cB$ the DG algebra $\cB=\End(E).$ Then these
pseudo-functors are also equivalent to the Maurer-Cartan
pseudo-functor $\cM\cC(\cB):\dgart \to \bf{Gpd}$ (I, Section 5).
\end{prop}

\begin{theo} Let $E, E^\prime \in C(\cM)$ be such that the DG
algebras $\End(E)$ and $\End (E^\prime)$ are quasi-isomorphic. (For
example assume that $E$ and $E^\prime$ are homotopy equivalent).
Then the pseudo-functors $\Def ^{\h} (E)$ ($ \simeq \coDef
^{\h}(E)$), $\Def ^{\h} (E^\prime)$ ($\simeq \coDef
^{\h}(E^\prime)$) are equivalent.
\end{theo}

\begin{theo}\label{def=def^h} Assume that $\cM$ has enough projectives and let $P\in
C^-(\cM)$ be h-projective such that $\Ext ^{-1}(P,P)=0.$ Then there
is an equivalence of deformation pseudo-functors
$$\Def _-(P)\simeq \Def ^{\h}_-(P).$$
\end{theo}

\begin{theo}\label{codef=codef^h}  Assume that $\cM$ has enough injectives and let $I\in
C^+(\cM)$ be h-injective such that $\Ext ^{-1}(I,I)=0.$ Then there
is an equivalence of deformation pseudo-functors
$$\coDef _-(I)\simeq \coDef ^{\h}_-(I).$$
\end{theo}

\begin{cor}\label{depends_on_RHom} Assume that $\cM$ has enough projectives (resp. injectives)
and $E\in D^-(\cM)$ (resp. $E\in D^+(\cM)$) is such that $\Ext
^{-1}(E,E)=0.$ Then the deformation pseudo-functor $\Def _-(E)$
(resp. $\coDef _- (E)$) depends up to an equivalence only on the
quasi-isomorphism class of the DG algebra $\bR \Hom (E,E).$ In
particular, let $\cN$ be another abelian category with enough
projectives (resp. enough injectives); assume that $F:\cM \to \cN$
is a functor such that $C(F):C(\cM)\to C(\cN)$ preserves
h-projectives (resp. h-injectives)  and induces an equivalence
$\tilde{F}:D^-(\cM)\to D^-(\cN)$ (resp. $\tilde{F}:D^+(\cM)\to
D^+(\cN)$),  then the pseudo-functors $\Def _-(E)$ and $\Def
_-(F(E))$ (resp. $\coDef _-(E)$ and $\coDef _-(F(E))$) are
equivalent.
\end{cor}

\begin{theo}\label{def=codef} Assume that $\cM$ has enough projectives and enough
injectives. Let $E\in D^b(\cM)$ be a complex with $\Ext
^{-1}(E,E)=0.$ Then there exists an equivalence of pseudo-functors
$$\Def _-(E)\simeq \coDef _-(E).$$
\end{theo}

These propositions and theorems will be proved in the next two
subsections.

\subsection{Relation with the homotopy deformation theory of DG modules}

Let $\cR$ be an artinian DG algebra. Recall (Remark
\ref{abelian_vs_dg}) that an $\cR ^{op}$-complex $M\in
C^{dg}_{\cR}(\cM)$ defines a DG functor from the DG category $\cM$
to the DG category of DG $\cR^{op}$-modules by the formula $N\mapsto
\Hom _{C^{dg}(\cM)}(N,M).$ This  may be considered also as a DG
module over the DG category $\cM ^{op}_{\cR}=(\cM \otimes
\cR)^{op}.$ Thus we obtain a full and faithful (by Yoneda) DG
functor
$$h^\bullet _{\cR}:C^{\dg}_{\cR}(\cM)\hookrightarrow \cM
^{op}_{\cR}\text{-mod}.$$

Notice that this DG functor $h^\bullet$ does not commute with the
DG functor $\phi ^*$ in general (because the usual Yoneda functor
$\cM \to \cM^{op}\text{-mod}$ is not exact). However we have the
following result which suffices for our purposes.

\begin{lemma}\label{commut} Let $\phi :\cR \to \cQ$ be a homomorphism of artinian
DG algebras.

a) There is a natural isomorphism of DG functors from $\cC
^{dg}_{\cQ}(\cM)$ to $\cM _{\cR}^{op}$-mod
$$h^\bullet _{\cR}\cdot \phi _*\simeq \phi _*\cdot h^\bullet_{\cQ}.$$

b) There is a natural isomorphism of DG functors from $\cC
^{dg}_{\cR}(\cM)$ to $\cM _{\cQ}^{op}$-mod
$$h^\bullet _{\cQ}\cdot \phi ^!\simeq \phi ^!\cdot h^\bullet _{\cR}.$$

c) There is a natural morphism of DG functors from $\cC
^{dg}_{\cR}(\cM)$ to $\cM _{\cQ}^{op}$-mod
$$\phi ^*\cdot h^\bullet _{\cR} \to h^\bullet _{\cQ}\cdot \phi ^*,$$
which is an isomorphism on objects $T\in C ^{dg}_{\cR}(\cM)$ such
that $T$ as a graded $\cR ^{op}$-module is isomorphic to $S\otimes
_kV,$ where $S\in C^{dg}(\cM)$ and $V$ is a finite dimensional $\cR
^{op}$-module. In particular it is an isomorphism on objects $T$
which are graded $\cR$-free.
\end{lemma}

\begin{proof} a) is obvious. Let us prove b) and c). Fix $M \in \cC
^{dg}_{\cR}(\cM)$ and $N\in \cM.$ We have
$$\phi ^! h ^M_{\cR}(N)=\Hom^{\cdot} _{\cR ^{op}}(Q,\Hom^{\cdot}
_{C^{dg}(\cM)}(N,M)),$$ and
$$ h ^{\phi ^!M}_{\cR}(N)=\Hom^{\cdot} _{C^{dg}(\cM)}(N,\Hom^{\cdot}
_{\cR ^{op}}(Q,M)).$$ Notice that both these complexes are
naturally isomorphic to
$$\Hom^{\cdot} _{C^{dg}(\cM)\otimes\cR ^{op}}(N\otimes Q,M).$$
This proves b).

Now
$$\phi ^* h _{\cR}^M(N)=\Hom^{\cdot}
_{C^{dg}(\cM)}(N,M)\otimes _{\cR}Q,$$ and
$$h _{\cR}^{\phi ^*M}(N)=\Hom^{\cdot} _{C^{dg}(\cM)}(N,M\otimes _{\cR}Q).$$
We define the morphism of complexes $\delta :\phi ^* h_{\cR}^M(N)\to
h _{\cR}^{\phi ^*M}(N)$ by
$$\delta (f\otimes q)(n)=(-1)^{\bar{q}\bar{n}}f(n)\otimes q.$$
Assume that $M$ as a graded object is isomorphic to $S\otimes V,$
where $S\in C^{dg}(\cM)$ and $V$ is a finite dimensional $\cR
^{op}$-module. Then both complexes $\phi ^* h _{\cR}^M(N)$ and $h
_{\cR}^{\phi ^*M}(N)$ are graded isomorphic to $\Hom^{\cdot}
_{C^{dg}(\cM)}(N,S)\otimes (V\otimes _{\cR} Q)$ and $\delta$ is the
identity map. This proves c) and the lemma.
\end{proof}

\noindent{\bf Warning.} In what follows we will compare
deformations and co-deformations of objects $E$ in $C(\cM)$ (or,
which is the same, in $C^{dg}(\cM)$) as defined above, with
deformations and co-deformations of DG modules over DG categories
(such as $h^E$ for example), as defined in \cite{ELOI}. These
pseudo-functors are denoted by the same symbols (like $\Def
^{\h}$), but we hope that there is no danger of confusion because
we always mention the corresponding argument (such as $E$ or
$h^E$).

\begin{cor}\label{compar_def^h} Fix $E\in C^{dg}(\cM).$  The collection of DG functors $\{h^\bullet _{\cR}\}$
defines morphisms of pseudo-functors
$$h^\bullet :\Def ^{\h}(E)\to \Def ^{\h} (h^E), \quad h^\bullet :\coDef ^{\h}(E)\to \coDef ^{\h}
(h^E)$$
\end{cor}

\begin{proof} Notice that for an artinian DG algebra $\cR$ the
graded $\cM ^{op}_{\cR}$-modules $h^E\otimes \cR$ and $h^{E\otimes
\cR}$ (resp. $h^E\otimes \cR ^*$ and $h^{E\otimes \cR ^*}$) are
naturally isomorphic. The rest follows from Lemma \ref{commut}.
\end{proof}

\begin{prop}\label{equiv_def^h} For every $E\in C^{dg}(\cM)$ the morphisms
$$h^\bullet :\Def ^{\h}(E)\to \Def ^{\h} (h^E), \quad h^\bullet :\coDef ^{\h}(E)\to \coDef
^{\h} (h^E)$$ are equivalences of pseudo-functors.
\end{prop}

\begin{proof} Since the DG functors $h^\bullet _{\cR}$ are full and faithful
it follows that the induced functor $h_{\cR}:\Def ^{\h}_{\cR}(E)\to
\Def ^{\h}_{\cR}(h^E)$ is full and faithful. It remains to show that
$h_{\cR}$ is essentially surjective.

Let $(\tilde{S},\id)\in \Def _{\cR}^{\h}(h^E).$ Consider the DG $\cM
^{op}$-module $p_*\tilde{S}.$ Notice that $\tilde{S}$ is just the DG
$\cM ^{op}$-module $p_*\tilde{S}$ together with a homomorphism of DG
algebras $\cR \to \End(p_*\tilde{S}).$ Thus it suffices to show that
$p_*\tilde{S}$ is isomorphic to $h^S$ for some $S\in C^{dg}(\cM).$
Notice that $p_*\tilde{S}$ is obtained from $E$ by taking finite
direct sums and iterated cone constructions. The DG functor
$h^\bullet$ preserves cones of morphisms, hence $p_*\tilde{S}$ is in
the essential image of the DG functor $h^\bullet.$

The same proof works for the pseudo-functors $\coDef ^{\h}.$
\end{proof}

\subsection{Proof of main properties of deformation pseudo-functors}

\begin{cor}\label{def^h=codef^h1} For any $E\in C(\cM)$ the pseudo-functors $\Def
^{\h}(E)$ and $\coDef ^{\h}(E)$ from $\dgart$ to ${\bf Gpd}$ are
equivalent.
\end{cor}

\begin{proof} Indeed, by the last proposition we have equivalences
$\Def ^{\h}(E)\simeq \Def ^{\h}(h^E)$ and $\coDef ^{\h}(E)\simeq
\coDef ^{\h}(h^E).$ It remains to apply I, Proposition 4.7.
\end{proof}

For $E\in C^{dg}(\cM)$ denote by $\cB$ the DG algebra $\End (E).$
Recall the Maurer-Cartan pseudo-functor $\cM\cC(\cB):\dgart \to {\bf
Gpd}$ (I, Definition 5.4).

\begin{cor}\label{def^h=MC} The pseudo-functors $\Def ^{\h}(E)$ ($\simeq \coDef
^{\h}(E)$) and $\cM\cC(\cB)$ are equivalent. In particular these
pseudo-functors depend (up to an equivalence) only on the
isomorphism class of the DG algebra $\End(E).$
\end{cor}

\begin{proof} This follows from Corollary \ref{def^h=codef^h1} and I, Proposition 6.1.
\end{proof}

Recall that for quasi-isomorphic DG algebras $\cB $ and $\cC$ the
corresponding Maurer-Cartan pseudo-functors $\cM\cC(\cB)$ and
$\cM\cC(\cC)$ are equivalent (I, Theorem 8.1). Hence we obtain the
following corollary.

\begin{cor}\label{depends_on_End} Assume that for $E, E^\prime \in C^{dg}(\cM)$ the DG
algebras $\End (E)$ and $\End (E^\prime)$ are quasi-isomorphic.
Then the pseudo-functors $\Def ^{\h}(E)$ ($\simeq \coDef
^{\h}(E)$)  and $\Def ^{\h}(E^\prime)$  ($\simeq \coDef
^{\h}(E^\prime)$) are equivalent.
\end{cor}

The next example is a copy of I, Proposition 8.3 in our context.

\begin{example}\label{when_End_quasi-iso} \rm{a) Assume that for $E, E^\prime \in C^{dg}(\cM)$ are
homotopy equivalent. Then the DG algebras $\End(E)$ and $\End
(E^\prime)$ are canonically quasi-isomorphic.

b) Let $P\in C^{dg}(\cM)$ and $I\in C^{dg}(\cM)$ be h-projective and
h-injective respectively. Assume that $f:P\to I$ is a
quasi-isomorphism. Then the DG algebras $\End(P)$ and $\End(I)$ are
canonically quasi-isomorphic.

The proof is the same as that of I, Proposition 8.3.}
\end{example}

We will need a more precise result as in Proposition 8.5 in Part I.

\begin{lemma}\label{homot_equiv} Fix an artinian DG algebra $\cR.$

 Let $g:E\to E^\prime $ be a homotopy equivalence in $C^{dg}(\cM).$
Assume that $(V,\id)\in \Def ^{\h}_{\cR}(E)$ and $(V^\prime,\id)\in
\Def ^{\h}_{\cR}(E^\prime)$ are objects that correspond to each
other via the equivalence $\Def ^{\h}_{\cR}(E)\simeq \Def
^{\h}_{\cR}(E^\prime)$ of Corollary \ref{depends_on_End}and Example
\ref{when_End_quasi-iso} a). Then there exists a homotopy
equivalence $\tilde{g}:V\to V^\prime$ which extends $g,$ i.e.
$i^*\tilde{g}=g.$ Similarly for the objects of $\coDef ^{\h}_{\cR}$
with $i^!$ instead of $i^*.$
\end{lemma}

\begin{proof} The full and faithful Yoneda DG functor $h^\bullet$
allows us to translate the problem to DG modules over the DG
category $\cM$ (Proposition \ref{equiv_def^h}). So it remains to
apply I, Proposition 8.5 a).
\end{proof}

\begin{theo}\label{def=def^h1} a) Assume that $\cM$ has enough projectives. Let $E\in
C^-(\cM)$ be such that $\Ext ^{-1}(E,E)=0.$ Choose an h-projective
$P\in C^-(\cM)$ and a quasi-isomorphism $P\to E.$ Then there exists
an equivalence of pseudo-functors
$$\Def _-(E)\simeq \Def ^{\h}_-(P).$$

b) Assume that $\cM$ has enough injectives. Let $E\in C^+(\cM)$ be
such that $\Ext ^{-1}(E,E)=0.$ Choose an h-injective $I\in C^+(\cM)$
and a quasi-isomorphism $E\to I.$ Then there exists an equivalence
of pseudo-functors
$$\coDef _-(E)\simeq \coDef ^{\h}_-(I).$$
\end{theo}

\begin{proof} We may and will assume that each $P^j\in \cM$ (resp.
$I^j\in \cM$) is projective (resp. injective).

We need a few preliminaries.

\begin{lemma}\label{h-proj_inj1} Fix $\cR\in \dgart _-.$ In the notation of the above theorem let
$(S,\sigma)\in \Def ^{\h}_{\cR}(P)$ (resp. $(S,\sigma)\in \coDef
^{\h}_{\cR}(I)$). Then $S\in C^{dg}_{\cR}(\cM)$ is h-projective
(resp. h-injective).
\end{lemma}

\begin{proof} Let $(S,\sigma)\in \Def ^{\h}_{\cR}(P).$ We may and
will assume that $i^*S=P$ and $\sigma =\id.$ By definition
$S^{\gr}\simeq (P\otimes \cR)^{\gr}$ and since $\cR$ is non-positive
the graded $\cR$-submodule $\oplus _{j>j_0}P^j\otimes \cR$ is
actually a subcomplex for each $j_0.$ Notice that for each $j$ the
$\cR^{op}$-complex $P^j\otimes \cR$ is h-projective (since $P^j$ is
projective). Hence also each $\cR ^{op}$-submodule $\oplus
_{j>j_0}P^j\otimes \cR \subset S$ is h-projective. Now we repeat the
argument in the proof of Proposition \ref{enough} a) to show that
$S$ is h-projective.

The proof for $(S,\sigma )\in \coDef ^{\h}_{\cR}(I)$ is similar.
\end{proof}

\begin{lemma}\label{lifting}  Let $\cR$ be an artinian
DG algebra and $S,T\in C^{dg}_{\cR}(\cM)$ be graded $\cR$-free
(resp. graded $\cR$-cofree).

a) There is an isomorphism of graded vector spaces
$\Hom^{\cdot}(S,T)=\Hom^{\cdot}(i^*S,i^*T) \otimes \cR$ (resp.
$\Hom^{\cdot}(S,T)=\Hom^{\cdot}(i^!S,i^!T) \otimes \cR$), which is
an isomorphism of graded algebras if $S=T.$ In particular, the map
$i ^*:\Hom^{\cdot} (S,T)\to \Hom^{\cdot} (i^*S,i^*T)$ (resp. $i
^!:\Hom^{\cdot}(S,T)\to \Hom^{\cdot}(i^!S,i^!T)$) is surjective.

b) The $\cR^{op}$-complex $S$ has a finite filtration with
subquotients isomorphic to $i^*S$ as objects in $C^{dg}(\cM)$
(resp. to $i^!S$ as objects in $C^{dg}_{\cR}(\cM)$).

c) The DG algebra $\End(S)$ has a finite filtration by DG ideals
with subquotients isomorphic to $\End (i^*S).$

d) If $f\in \Hom^{\cdot} (S,T)$ is a closed morphism of degree
zero such that $i^*f$ (resp. $i^!f$) is an isomorphism or a
homotopy equivalence, then $f$ is also such.
\end{lemma}

\begin{proof} The full and faithful Yoneda DG functor $h^\bullet$
allows us to translate the problem to DG modules over the DG
category $\cM$ (Lemma \ref{commut}). So it remains to apply I,
Proposition 3.12.
\end{proof}

Now we can prove the theorem. We first prove a). Using Remark
\ref{quasi-iso_def} it suffices to prove that the pseudo-functors
$\Def _-(P)$ and $\Def _-^{\h}(P)$ are equivalent. Let us define a
morphism
$$\beta : \Def _-^{\h}(P)\to \Def _-(P)$$
Fix $\cR \in \dgart _-$ and let $(S,\sigma )\in \Def
_{\cR}^{\h}(P).$ By Lemma \ref{h-proj_inj1} the $\cR^{op}$-complex
$S$ is h-projective. Hence $\bL i^*S=i^*S$ and therefore $(S,\sigma
)\in \Def _{\cR}(P).$ This defines a functor $\beta _{\cR}:\Def
_{\cR}^{\h}(P)\to \Def _{\cR}(P)$ and a morphism of pseudo-functors
$\beta : \Def _-^{\h}(P)\to \Def _-(P).$ We need to show that $\beta
_{\cR}$ is an equivalence.

\medskip

\noindent{\bf Surjective on isomorphism classes.} Let $(T,\tau )\in
\Def _{\cR}(P).$ We may and will assume that $T\in
C^{dg}_{\cR}(\cM)$ is h-projective and graded $\cR$-free
(Proposition \ref{enough} a)). Thus $(T,\tau )\in \Def
_{\cR}^{\h}(i^*T).$ Since $T$ is h-projective, so is $i^*T$ and
hence $\tau :\bL i^*T=i^*T\to P$ is a homotopy equivalence. It
follows from Lemma \ref{homot_equiv} that there exists $(S,\id)\in
\Def _{\cR}^{\h}(P)$ such that $S$ and $T$ are homotopy equivalent
and $(S,\id)$ and $(T,\tau)$ are isomorphic objects in $\Def
_{\cR}(P).$ I. e. $\beta _{\cR}(S,\id)\simeq (T,\tau).$

\medskip

\noindent{\bf Full.} Let $(S,\id ), (S^\prime, \id )\in \Def
^{\h}_{\cR}(P).$ Let $f:\beta_{\cR}(S,\id ) \to
\beta_{\cR}(S^\prime, \id )$ be an isomorphism in $\Def _{\cR}(P).$
Since $S,S^\prime$ are h-projective (Lemma \ref{h-proj_inj1}) this
isomorphism $f$ is a homotopy equivalence. Because $P$ is
h-projective $i^*f$ is homotopic to $\id _P.$ Let $h:i^*f\to \id$ be
a homotopy. Since $S,$ $S^\prime$ are graded $\cR$-free the map
$i^*:\Hom^{\cdot} (S, S^\prime )\to \Hom^{\cdot} (P,P)$ is
surjective (Lemma \ref{lifting} a)). Choose a lift $\tilde{h}:S\to
S^\prime[1]$ of $h$ and replace $f$ by $\tilde{f}=f-d\tilde{h}.$
Then $i^*\tilde{f}=id.$ Since $S$ and $S^\prime$ are graded
$\cR$-free $\tilde{f}$ is an isomorphism (Lemma \ref{lifting} d)).
Thus $\tilde{f}:(S,\id ) \to (S^\prime, \id )$ is a morphism in
$\Def _{\cR}^{\h}(P)$ such that $\beta _{\cR}\tilde{f}=f.$

\medskip

\noindent{\bf Faithful.} Let $(S,\id ),(S^\prime, \id )\in \Def
_{\cR}^{\h}(P)$ and let $g_1,g_2:S\to S^\prime$ be two isomorphisms
(in $C^{dg}_{\cR}(\cM)$) such that $i^*g_1=i^*g_2=\id _P.$ That is
$g_1,g_2$ are maps which represent morphisms in $\Def
_{\cR}^{\h}(P).$ Assume that $\beta _{\cR} (g_1)=\beta _{\cR}
(g_2),$ i.e. there exists a homotopy $s:g_1\to g_2.$ Then
$d(i^*s)=i^*(ds)=0.$ Since by our assumption $H^{-1}\Hom^{\cdot}
(P,P)=0$ there exists $t\in \Hom ^{-2}(P,P)$ with $dt=i^*s.$ Choose
a lift $\tilde{t}\in \Hom ^{-2}(S,S^\prime)$ of $t.$ Then
$\tilde{s}:=s-d\tilde{t}$ is an allowable homotopy between $g_1$ and
$g_2.$ This proves that $\beta _{\cR} $ is faithful.

The proof of part b) of the theorem is similar and we omit it.
\end{proof}

\begin{theo}\label{def=codef1} Assume that $\cM$ has enough projectives and enough
injectives. Let $E\in D^b(\cM)$ be a complex such that $\Ext
^{-1}(E,E)=0.$ Then there exists an equivalence of pseudo-functors
$$\Def _-(E)\simeq \coDef _-(E).$$
\end{theo}

\begin{proof} Choose quasi-isomorphisms $P\to E$ and $E\to I,$ where
$P$ is a bounded above h-projective and $I$ is a bounded below
h-injective (Proposition \ref{enough}). Then by Theorem
\ref{def=def^h1} there exist equivalences of pseudo-functors
$$\Def _-(E)\simeq \Def ^{\h}_-(P),\quad \coDef _-(E)\simeq \coDef
_-^{\h}(I).$$ But pseudo-functors $\Def ^{\h}_-(P)$ and $\coDef
_-^{\h}(I)$ are equivalent by Example \ref{when_End_quasi-iso},
Corollary \ref{depends_on_End}.
\end{proof}

\begin{cor}\label{depends_on_RHom1}  Assume that $\cM$ has enough projectives (resp. injectives) and $E\in
D^-(\cM)$ (resp. $E\in D^+(\cM)$) is such that $\Ext ^{-1}(E,E)=0.$
Then the deformation pseudo-functor $\Def _-(E)$ (resp. $\coDef _-
(E)$) depends up to an equivalence only on the quasi-isomorphism
class of the DG algebra $\bR \Hom (E,E).$ In particular, let $\cN$
be another abelian category with enough projectives (resp. enough
injectives); assume that
 $F:\cM \to \cN$ is a functor such that $C(F):C(\cM)\to C(\cN)$ preserves
h-projectives (resp. h-injectives)  and induces an equivalence
$\tilde{F}:D^-(\cM)\to D^-(\cN)$ (resp. $\tilde{F}:D^+(\cM)\to
D^+(\cN)$),  then the pseudo-functors $\Def _-(E)$ and $\Def
_-(F(E))$ (resp. $\coDef _-(E)$ and $\coDef _-(F(E))$) are
equivalent.
\end{cor}

\begin{proof}  This follows from Theorem \ref{def=def^h1} and Corollary \ref{depends_on_End}.
\end{proof}

\subsection{Relation with the derived deformation theory of DG modules} Finally in the situation of
Theorem \ref{def=def^h1} we want to interpret the {\it derived}
deformation pseudo-functors $\Def _- $ and $\coDef _-$ in a
context of DG modules.

\begin{theo}\label{compar_def} Assume that $\cM$ has  enough projectives and let $P\in
C^-(\cM)$ be h-projective such that $\Ext ^{-1}(P,P)=0.$ Then the DG
functor $h^\bullet:C^{dg}(\cM)\to \cM ^{op}\text{-mod}$ induces an
equivalence of pseudo-functors $\Def _-(P)\simeq \Def _-(h^P).$
\end{theo}

\begin{proof} By Proposition \ref{equiv_def^h} the morphism of pseudo-functors
$$h^\bullet :\Def ^{\h}_-(P)\to \Def ^{\h}_-(h^P)$$
is an equivalence. By Theorem \ref{def=def^h1} a) there exists an
equivalence of pseudo-functors $\Def _-(P)\simeq \Def ^{h}_-(P).$ We
claim that the pseudo-functors $\Def _-(h^P)$ and $\Def
^{\h}_-(h^P)$ are also equivalent. Indeed, notice that the DG $\cM
^{op}$-module $h^P$ satisfies property (P) (Definition 3.2 in Part
I). Hence it is h-projective. Therefore
$$\Ext ^{-1}(h^P,h^P)=H^{-1}\Hom^{\cdot} (h^P,h^P)\simeq H^{-1}\Hom^{\cdot} (P,P)=0.$$
Clearly $h^P$ is bounded above. Hence by I, Theorem 11.6 a) we have
$\Def _-(h^P)\simeq \Def ^{\h}_-(h^P).$

Combining these three equivalences we obtain $\Def _-(P)\simeq \Def
_-(h^P).$
\end{proof}

\begin{remark} \rm{Notice that the DG functor $h^\bullet :C^{dg}(\cM)\to
\cM ^{op}\text{-mod}$ does not preserve quasi-isomorphisms in
general. If $\cM$ has enough projectives then we can consider a
similar DG functor
$${}^Ph^\bullet :C^{dg}(\cM)\to \cP ^{op}\text{-mod},$$
where $\cP \subset \cM$ is the full subcategory of projectives.
This DG functor ${}^Ph^\bullet $ has all the good properties of
$h^\bullet$ (full and faithful, induces an equivalence of homotopy
deformation and co-deformation pseudo-functors, etc.) and in
addition it preserves quasi-isomorphisms. Thus ${}^Ph^\bullet$ is
better suited than $h^\bullet$ for comparing derived deformation
pseudo-functors.}
\end{remark}

Next we want to prove the analogue of Theorem \ref{compar_def} for
the pseudo-functor $\coDef _-$ in case $\cM $ has enough
injectives. We can only prove it with an extra finiteness
assumption.

We are still going to work with a covariant DG functor from
$C^{gd}(\cM)$ to $\cM ^{op}\text{-mod},$ but it will not be
$h^\bullet.$ Consider the DG functor $h_\bullet ^*:C^{dg}(\cM)\to
\cM ^{op}\text{-mod}$ defined by
$$h^*_M(N):=\Hom^{\cdot} _{C^{dg}(\cM)}(M,N)^*,$$
where $(\cdot )*$ denotes the (graded) $k$-dual. Recall (I, Section
3.1) that for any $W\in \cM ^{op}\text{-mod},$ $M\in C^{dg}(\cM)$
$$\Hom^{\cdot} (W,h^*_M)=W(M)^*.$$
In particular
$$\Hom^{\cdot} (h^*_{M_1},h^*_{M_2})=h^*_{M_1}(M_2)^*=\Hom^{\cdot}
_{C^{dg}(\cM)}(M_1,M_2)^{**}.$$ Therefore the DG functor
$h^*_\bullet$ is not full in general, but it induces a
quasi-isomorphism
$$h^*_\bullet :\Hom^{\cdot} (M_1,M_2)\to \Hom^{\cdot} (h^*_{M_1},h^*_{M_2})$$
if $\dim H^i\Hom (M_1,M_2)<\infty$ for all $i.$

Also for each $M\in \cM$ the DG $\cM ^{op}$-module $h_M^*$ is
h-injective. Hence $h^*_M$ is h-injective for each $M\in C^+(\cM).$

\begin{theo}\label{compar_codef} Assume that $\cM$ has enough injectives and let $I\in
C^+(\cM)$ be h-injective such that $\Ext ^{-1}(I,I)=0.$ Assume that
for each $i$ $\dim \Ext ^i(I,I)<\infty.$ Then the DG functor
$$h_\bullet ^*:C^{dg}(\cM)\to \cM ^{op}\text{-mod}$$
induces an equivalence of pseudo-functors $\coDef _-(I)\simeq \coDef
_-(h^*_I).$
\end{theo}

\begin{proof} By Theorem \ref{def=def^h1} b) there exists an equivalence of
pseudo-functors $$\coDef _-(I)\simeq \coDef ^{\h}_-(I).$$

Since $\dim H^i\Hom (I,I)<\infty $ for each $i$ the homomorphism
of DG algebras $h_\bullet ^*:\End (I)\to \End (h^*_I)$ is a
quasi-isomorphism. Hence the pseudo-functors $\coDef ^{\h}_-(I)$
and $ \coDef ^{\h}_-(h^*_I)$ are equivalent (Corollary
\ref{def^h=MC} and I, Proposition 6.1, I, Theorem 8.1).

The DG $\cM ^{op}$-module $h_I^*$ is h-injective and bounded
below. Hence by I, Theorem 11.6 b)
$$\coDef _-^{\h}(h^*_I)\simeq \coDef _-(h^*_I).$$
Combining the above three equivalences we obtain the desired
equivalence
$$\coDef _-(I)\simeq \coDef _-(h_I^*).$$
\end{proof}

In case of finite injective dimension we could still use the DG
functor $h^\bullet $ to compare the derived co-deformation
pseudo-functors. Namely, we have the following result.

\begin{theo}\label{compar_codef1} Assume that $\cM$ has enough injectives and let $I\in
C^b(\cM)$ be a (bounded) h-injective such that $\Ext ^{-1}(I,I)=0.$
Then the DG functor
$$h^\bullet :C_{dg}(\cM)\to \cM^{op}\text{-mod}$$
induces an equivalence of pseudo-functors
$$\coDef _-(I)\simeq \coDef _-(h^I).$$
\end{theo}

\begin{proof} By Proposition \ref{equiv_def^h} we have an equivalence of
pseudo-functors
$$h^\bullet :\coDef ^{\h}(I)\to \coDef ^{\h}(h^I).$$

By Theorem \ref{codef=codef^h} there is an equivalence of
pseudo-functors
$$\coDef _-(I)\simeq \coDef ^{\h}_-(I).$$

Finally, notice that $h^I$ is a bounded h-projective object in $\cM
^{op}\text{-mod}.$ Hence by I, Theorem 11.6 b) there is an
equivalence of pseudo-functors
$$\coDef ^{\h}_-(h^I)\simeq \coDef _-(h^I).$$
This proves the theorem.
\end{proof}

\part{Geometric examples, applications and conjectures}

\section{Deformations of bounded complexes on locally noetherian
schemes}

Fix a locally noetherian scheme $X$ over $k.$ Let $E$ be a bounded
complex of quasi-coherent sheaves on $X.$ The abelian category
$\Qcoh _X$ of quasi-coherent sheaves on $X$ does not have enough
projectives in general. Still there is a natural (classical) derived
deformation pseudo-functor
$$\Def _{\cl}(E):\cart \to {\bf Gpd},$$
which is defined using h-flat objects.

The abelian category $\Qcoh _X$ has enough injectives, so we can
define the pseudo-functor $\coDef_-(E)$ as in Definition
\ref{codef_-} above. Our main result (Theorem \ref{def_cl=codef_cl}
below) claims that the pseudo-functors $\Def _{\cl}(E)$ and $\coDef
_{\cl}(E)$ are naturally equivalent. This allows us to consider the
pseudo-functor $\coDef _-(E)$  as a natural extension to $\dgart _-$
of the classical deformation functor $\Def _{\cl}(E).$

Let us first introduce some notation. For a scheme $Y$ we denote by
$\Mod _Y,$ $\Qcoh _Y,$ $\coh _Y$ the abelian categories of $\cO
_Y$-modules, quasi-coherent $\cO _Y$-modules and coherent $\cO
_Y$-modules respectively. Denote by $D(Y),$ $D(\Qcoh _Y),$ $D(\coh
_Y)$ the corresponding derived categories and by $D^{\pm}_{\Qcoh
}(Y),$ $D^b_{\coh}(Y),$ $D^{\pm}_{\coh} (\Qcoh _Y),$ ... their usual
full subcategories defined by a cohomological condition.

Note that none of the categories $\Mod _Y,$ $\Qcoh _Y,$ $\coh _Y$
has enough projectives in general. The categories $\Mod _Y,$ $\Qcoh
_Y$ have enough injectives and if the scheme $Y$ is locally
noetherian, then the natural functor
$$D^+(\Qcoh _Y)\to D^+_{\Qcoh}(Y)$$
is an equivalence of categories \cite{Hart}.

\begin{defi}\label{h-flat} A complex $F\in C(\Mod _Y)$ is h-flat if the complex
$F\otimes _{\cO _Y}G$ is acyclic for every acyclic $G\in C(\Mod
_Y).$
\end{defi}

For every $S\in C(\Mod _Y)$ Spantelstein in \cite{Sp} has
constructed a functorial h-flat resolution. That is he defines an
h-flat $F(S)\in C(\Mod _Y)$ and a quasi-isomorphism $F(S)\to S.$ The
complex $F(S)$ consists of $\cO _Y$-modules which are direct sums of
sheaves $\cO _U$ for affine open subsets $U\subset Y$ ($\cO _U$ is
the extension by zero to $Y$ of the structure sheaf of $U$). Using
these h-flat resolutions we may define derived functors $\bL
f^*:D(Y)\to D(Z)$ for a morphism of schemes $f:Z\to Y.$ Namely, put
$$\bL f^*(S):=f^*(F(S)).$$

For a commutative local artinian algebra $\cR$ and a scheme $Y$ put
$$Y_{\cR}=X\otimes _{\Spec k}\Spec \cR$$
and denote by
$i:Y\hookrightarrow Y_{\cR}$ the closed embedding.

\begin{defi}\label{def_cl} Let $X$ be a scheme, $E\in D_{\Qcoh}(X).$ We
define the pseudo-functor
$$\Def _{\cl}(E):\cart \to {\bf Gpd}$$
of "classical" deformations of $E$ as follows.

Fix a commutative local artinian algebra $\cR.$ An object of the
groupoid $\Def _{\cR}(E)$ is a pair $(S,\sigma),$ where $S\in
D_{\Qcoh}(X_{\cR})$ and $\sigma :\bL i^*S\to E$ is an isomorphism in
$D_{\Qcoh}(X).$ A morphism between two such pairs $(S,\sigma )$ and
$(S^\prime ,\sigma ^\prime)$ is an isomorphism $f:S\to S^\prime$
such that $\sigma =\sigma ^\prime \cdot \bL i^*(f).$

A homomorphism $\phi :\cR \to \cQ$ of commutative local artinian
algebras induces a morphism of schemes $\phi :X_Q\to X_R$ which fits
in a commutative diagram
$$\begin{array}{rcl} X_Q & \stackrel{\phi}{\to} & X_R\\
i\uparrow & & \uparrow i\\
X & = & X.\end{array}
$$
Hence we obtain the functor $\bL \phi ^*:\Def _{\cR}(E)\to \Def
_{\cQ}(E).$ This defines the pseudo-functor $\Def _{\cl}(E):\cart
\to {\bf Gpd}.$
\end{defi}

In \cite{Sp} it is also shown that for every object $S\in C(\Mod
_Y)$ there exists an h-injective $J\in C(\Mod _Y)$ and a
quasi-isomorphism $S\to J.$ As usual we define right derived
functors using h-injectives.

For example if $\phi :\cR \to \cQ$ is a homomorphism of commutative
artinian local algebras we obtain the functor $$\phi ^!:C(\Mod
_{Y_{\cR}})\to C(\Mod _{Y_{\cQ}}),\quad \phi ^!T:=\cH om _{\cO
_{Y_{\cR}}}(\cO _{Y_{\cQ}}, T)$$ and its derived functor $\bR \phi
^!:D(Y_{\cR})\to D(Y_{\cQ}).$

In particular for a commutative local artinian algebra $\cR$ and a
scheme $Y$ the closed embedding $i:Y\hookrightarrow Y_{\cR}$ induces
the functor
$$i^!:C(\Mod _{Y_{\cR}})\to C(\Mod _Y), \quad i^!T:=\cH om_{\cO
_{Y_{\cR}}}(\cO _Y,T),$$ and its derived functor
$$\bR i^!:D(Y_{\cR})\to D(Y).$$

\begin{defi}\label{codef_cl} Let $X$ be a scheme (over $k$) and $E\in D_{\Qcoh}(X).$
We define the pseudo-functor
$$\coDef _{\cl}(E):\cart \to {\bf Gpd}$$
of ("classical") derived co-deformations of $E$ as follows.

Let $\cR$ be a commutative local artinian algebra. An object of the
groupoid $\coDef _{\cR}(E)$ is a pair $(T, \tau),$ where $T\in
D_{\Qcoh}(X_{\cR})$ and $\tau :E\to \bR i^!T$ is a
quasi-isomorphism. A morphism between two such object $(T,\tau )$
and $(T^\prime ,\tau ^\prime)$ is a quasi-isomorphism $f:T\to
T^\prime$ such that $\tau ^\prime =\bR i^!(f)\cdot \tau.$

A homomorphism $\phi :\cR \to \cQ$ of commutative local artinian
algebras induces the morphism $X_{\cQ}\to X_{\cR}$ and hence the
functor
$$\bR \phi ^!:\coDef _{\cR}(E)\to \coDef _{\cQ}(E).$$
thus we obtain the pseudo-functor $\coDef _{\cl}(E):\cart \to {\bf
Gpd}.$
\end{defi}

The next theorem is our main result of this section.

\begin{theo}\label{def_cl=codef_cl} Let $X$ be a locally noetherian scheme, $E\in
D^b_{\Qcoh}(X).$ Then there exists an equivalence of pseudo-functors
$$\Def _{\cl}(E)\to \coDef _{\cl}(E).$$
\end{theo}

This theorem follows from a more precise Theorem 4.16 below.

We need a few lemmas.

Fix a commutative artinian local algebra $\cR.$ Denote as usual by
$\cR ^*$ the $\cR$-module $\Hom^{\cdot} _k(\cR ,k).$ Let $F_{\cR}$
and $I_{\cR}$ denote the categories of free and injective
$\cR$-modules respectively.

\begin{lemma}\label{str_of_inj} a) $\cR ^*$ is the unique (up to isomorphism)
 indecomposable injective $\cR$ module.

b) Every injective $\cR$-module is isomorphic to a direct sum of
copies of $\cR ^*.$ A direct sum of injective $\cR$-modules is
injective.

c) The categories $F_{\cR}$ and $I_{\cR}$  are equivalent. The
mutually inverse equivalences are given by $\varphi :F_{\cR}\to
I_{\cR},$ $\psi :I_{\cR}\to F_{\cR},$ where
$$\varphi (M)=M\otimes _{\cR}\cR ^*,\quad \psi (N)=\Hom^{\cdot} _{\cR}(\cR
^*,N).$$

d) The functorial diagram
$$\begin{array}{ccc}
F_{\cR} & \stackrel{\varphi}{\to} & I_{\cR}\\
i^*\downarrow & & \downarrow i^!\\
k\text{-mod} & = & k\text{-mod}
\end{array}$$
commutes.

\end{lemma}

\begin{proof} a) Since the ring $\cR$ is noetherian and has a unique
prime ideal $m\subset \cR,$ the $\cR$-module $\cR ^*$ (which is the
injective hull of $\cR /m=k$) is the unique (up to isomorphism)
indecomposable injective $\cR$-module (\cite{Mat}).

b) This follows from a) and the fact that the abelian category of
$\cR$-modules is locally noetherian (\cite{Hart}, \cite{Gab}).

c) Notice that the natural map of $\cR$-modules $\cR \to
\Hom^{\cdot} _{\cR}(\cR ^* ,\cR ^*)$ is an isomorphism. Now
everything follows from b) and the fact that the functors $\varphi
$ and $\psi$ commute with direct sums.

d) Let $M$ be an $\cR$-module. Define a morphism of vector spaces
$$\beta :M\otimes _{\cR}k\to \Hom^{\cdot} _{\cR}(k,M\otimes _{\cR}\Hom^{\cdot}
_k(\cR ,k)), \quad \beta (m\otimes \xi)(\eta)=m\otimes \xi \eta
\epsilon,$$ where $\epsilon :\cR \to k$ is the augmentation map.
This map is an isomorphism if $M=\cR.$ Hence it is an isomorphism
for a free $\cR$-module $M.$
\end{proof}

\begin{defi}\label{R-free_sheaves} An $\cQ _{X_{\cR}}$-module $M$ is called $\cR$-free
(resp. $\cR$-injective) if every stalk $M_x$ is free (resp.
injective) as an $\cR$-module. We call a complex $S\in C(\Mod
_{X_{\cR}})$ $\cR$-free (resp. $\cR$-injective) if every $\cO
_{X_{\cR}}$-module $S^j$ is such. Denote by $C_F(\Mod _{X_{\cR}})$
(resp. $C_I(\Mod _{X_{\cR}})$) the full subcategories of $C(\Mod
_{X_{\cR}})$ which consist of $\cR$-free (resp. $\cR$-injective)
complexes.
\end{defi}

\begin{prop}\label{some_functors} Consider the functors $\varphi ,\psi :C(\Mod
_{X_{\cR}})\to C(\Mod _{X_{\cR}}).$
$$\varphi(S)=S\otimes _{\cO _{X_{\cR}}}p^!\cO _X=S\otimes _{\cR}\cR ^*,\quad \psi(T)=\cH om
_{\cO_{X_{\cR}}}(p^!\cO _X,T)=\Hom^{\cdot} _{\cR}(\cR ^*,T),$$ where
$p^!\cO _X=\cO _{X}\otimes _k\cR ^*.$

a) These functors induce mutually inverse equivalences of categories
$$\varphi :C_F(\Mod _{X_{\cR}}) \to C_I(\Mod _{X_{\cR}}),\quad
\psi:C_I(\Mod _{X_{\cR}})\to C_F(\Mod _{X_{\cR}}).$$

b) The functorial diagram
$$\begin{array}{rcl} C_F(\Mod _{X_{\cR}}) & \stackrel{\phi}{\to}
 & C_I(\Mod _{X_{\cR}}) \\
i^* \downarrow & & \downarrow i^!\\
C(\Mod _X) & =& C(\Mod _X)
\end{array}$$
commutes.
\end{prop}

\begin{proof} a) Let $x\in X_{\cR}.$ Then $\cO _{X_{\cR},x}=\cO
_{X,x}\otimes _k\cR.$ We have
$$\varphi (S)_x=S_x\otimes _{\cR}\cR ^*$$ and $$\psi(T)_x=\Hom^{\cdot}
_{\cR}(\cR ^* ,T_x).$$

Now the assertion follows from part c) of Lemma \ref{str_of_inj}.

b) For an $\cO _{X_{\cR}}$-module $M$ we have
$$i^*M=M\otimes _{\cO _{X_{\cR}}}\cO _X=M\otimes _{\cR}k,$$
$$i^!M=\cH om _{\cO _{X_{\cR}}}(\cO _X,M)=\Hom^{\cdot} _{\cR}(k,M).$$
Hence the assertion follows from part d) of Lemma
\ref{str_of_inj}.
\end{proof}

\begin{prop}\label{F---bounded} Let $F\in C(\Mod _{X_{\cR}}).$ Suppose
that $\bL i^*F\in D^b(X).$ Then $F$ is quasi-isomorphic to a bounded
$\cR$-free complex.
\end{prop}

Before we prove the proposition let us state an immediate corollary.

\begin{cor}\label{bounded_in_def} Given $(S,\sigma)\in \Def _{\cR}(E)$ there exists an
isomorphic $(S^\prime ,\sigma ^\prime ) \in \Def _{\cR}(E)$ such
that $S^\prime $ is a bounded $\cR$-free complex.
\end{cor}

\begin{proof} Let us prove the proposition. This is done in the next
two lemmas.

\begin{lemma}\label{h-flat_h-inj} For every $S\in C(\Mod _{X_{\cR}})$ there exist
quasi-isomorphisms $P\to S$ and $S\to J,$ where $P$ is h-flat and
$\cR$-free and $J$ is h-injective and $\cR$-injective.
\end{lemma}

\begin{proof} This is proved in \cite{Sp}. Namely, the assertion about
$P$ follows from Proposition 5.6 and that about $J$ follows from
Lemma 4.3 and Theorem 4.5 in \cite{Sp}.
\end{proof}

Let $F$ be as in the proposition. By the above lemma we may and will
assume that $F$ is h-flat and $\cR$-free. Hence $\bL i^*F=i^*F.$ The
following claim implies the proposition.

\medskip

\noindent{\bf Claim.} Let $K\in C(\Mod _{X_{\cR}})$ be $\cR$-free
and such that $H^j(i^*K)=0$ for $j<j_0$ and $j>j_1.$ Then $K$ is
quasi-isomorphic to its truncation $\tau _{\leq j_1}\tau _{\geq
j_0}F$ and moreover this truncation is $\cR$-free.

\medskip

Our claim follows from the next lemma.

\begin{lemma}\label{suffices_mod_m} Let
$M^\bullet:=M^{-1}\stackrel{d^{-1}}{\to}M^0\stackrel{d^0}{\to}M^1$
be a complex of free $\cR$-modules. Assume that
$H^0(M^\bullet\otimes _{\cR}k)=0.$ Then $H^0(M^\bullet)=0$ and $\Ker
d^0$ is a free $\cR$-module.
\end{lemma}

\begin{proof} We can find a finite filtration $\cR \supset m_1\supset
m_2...$ by ideals such that $m_s/m_{s+1}\simeq k.$ Consider the
induced filtration on the complex $M^\bullet$:
$$M^\bullet\supset m_1M^\bullet\supset m_2M^\bullet...$$
Then each subquotient $m_sM^\bullet/m_{s+1}M^\bullet$ is isomorphic
to the complex $M^\bullet \otimes _{\cR}k.$ Hence
$H^0(m_sM^\bullet/m_{s+1}M^\bullet)=0$ for each $s$ and hence
$H^0(M^\bullet)=0$ by devissage. This proves the first assertion of
the lemma.

To prove the second one we use the following fact:
 an $\cR$-module $N$ is free if (and only if) $\Tor
^{\cR}_1(N,k)=0$ (\cite{AM}, Ch. 2 Ex. 26 and \cite{AnRo} Prop.
2.1.4).

Consider the exact sequence
$$M^{-1}\to M^0\to M^1\to \coker d^0\to 0.$$
Then by our assumption $H^0(M ^\bullet \otimes _{\cR}k)=\Tor
^{\cR}_1(\coker d^0,k)=0.$ Hence $\coker d^0$ is a free
$\cR$-module. Thus $\im d^0$ is free and hence also $\Ker d^0$ is
free. This proves the lemma.
\end{proof}

The Proposition \ref{F---bounded} is proved.
\end{proof}

Now we want to prove the analogue of Corollary
\ref{bounded_in_def} for the co-deformation functor.

\begin{prop}\label{G---bounded} Let $G\in C(\Mod _{X_{\cR}}).$
Assume that $\bR i^!G\in D^b(X).$ Then $G$ is quasi-isomorphic to a
bounded $\cR$-injective complex.
\end{prop}

\begin{proof} By Lemma \ref{h-flat_h-inj} we may and will assume that $G$ is
h-injective and $\cR$-injective. Hence $\bR i^!G=i^!G.$ Then by
Proposition \ref{some_functors} the complex $\psi (G)$ is $\cR$-free
and $i^*\psi (G)=i^!G,$ so we may and will assume that $H^j(i^*\psi
(G))=0$ for $j<j_0$ and $j>j_1.$ By the Claim in the proof of
Proposition \ref{F---bounded} the complex $\psi (G)$ is
quasi-isomorphic to its truncation $\tau _{\leq j_1}\tau _{\geq
j_0}\psi (G)$ which is moreover $\cR$-free. But then this truncation
is a direct summand of $\psi (G)$ as a complex (of sheaves) of free
$\cR$-modules. Applying the functor $\varphi$ from Proposition
\ref{some_functors} we find that $G=\varphi (\psi (G))$ is
quasi-isomorphic to its truncation $\varphi (\tau _{\leq j_1}\tau
_{\geq j_0}\psi (G))=\tau _{\leq j_1}\tau _{\geq j_0}\varphi (\psi
(G))$ which is moreover $\cR$-injective. This proves the
proposition.
\end{proof}

We obtain the immediate corollary.

\begin{cor}\label{bounded_in_codef} Given $(T,\tau)\in \coDef _{\cR}(E)$ there exists an
isomorphic $(T^\prime ,\tau ^\prime ) \in \coDef _{\cR}(E)$ such
that $T^\prime $ is a bounded $\cR$-injective complex.
\end{cor}

\begin{prop}\label{F,G---acyclic} a) Let $F\in C^b(\Mod _{X_{\cR}})$ be a bounded $\cR$-free complex. Then $F$ is
acyclic for the functors $i^*$ and $\varphi.$ That is $i^*F=\bL
i^*F$ and $\varphi (F)=\bL \varphi (F).$

b) Let $G\in C^b( \Mod _{X_{\cR}})$ be a bounded $\cR$-injective
complex. Then $G$ is acyclic for the functors $i^!$ and $\psi.$ That
is $i^!G=\bR i^!G,$ $\psi (G)=\bR \psi (G).$

c) If $F$ (resp. $G$) has quasi-coherent cohomology then so do
$i^*F$ and $\varphi(F)$ (resp. $i^!G$ and $\psi (G)$).
\end{prop}

\begin{proof} a) We have $i^*F=F\otimes _{\cR}k$ and
$\varphi (F)=F\otimes _{\cR}\cR^*.$ Now use the fact that a bounded
complex of free $\cR$-modules is h-projective.

b) We have $i^!G=\Hom^{\cdot} _{\cR}(k,G)$ and $\psi
(G)=\Hom^{\cdot} _{\cR}(\cR ^*,G).$ Now use the fact that a bounded
complex of injective $\cR$-modules is h-injective.

c) This follows from a), b) and Propositions 3.3 and 4.3 in
\cite{Hart}.
\end{proof}

\begin{theo}\label{equivalences} Let $X$ be a locally noetherian scheme and $E\in
D^b_{\Qcoh}(X).$ Then there exist mutually inverse equivalences of
pseudo-functors
$$\bL \varphi :\Def _{\cl}(E)\to \coDef _{\cl}(E),$$
$$\bR \psi :\coDef _{\cl}(E)\to \Def _{\cl}(E),$$ such that for a
commutative  artinian local algebra $\cR$ and $S=(S,\sigma)\in \Def
_{\cR}(E),$ $T=(T,\tau )\in \coDef _{\cR}(E)$
$$\bL \phi (S)=S\stackrel{\bL}{\otimes }_{\cR}\cR ^*,\quad \bR \psi (T)=\bR \cH om _{\cR}(\cR ^*,T).$$
\end{theo}

\begin{proof} Fix an artinian commutative local algebra $\cR.$ By
Proposition \ref{F---bounded} (resp. Proposition \ref{G---bounded})
the category $\Def _{\cR}(E)$ (resp. $\coDef _{\cR}(E)$) is
equivalent to its full subcategory consisting of objects $(S,\sigma
)$ (resp. $(T,\tau )$) such that $S\in C^b(\Mod _{X_{\cR}})$ is
$\cR$-free (resp. $T\in C^b(\Mod _{X_{\cR}})$ is $\cR$-injective).
Moreover by Proposition \ref{F,G---acyclic} $\bL i^*S=i^*S,$ $\bR
i^!T=i^!T$ and $\bL \phi (S)=S\otimes _{\cR}\cR ^*,$ $\bR \psi
(T)=\Hom^{\cdot} _{\cR}(\cR ^* ,T).$ Now the theorem follows from
Proposition \ref{some_functors} and Proposition \ref{F,G---acyclic}
c). This also proves Theorem \ref{def_cl=codef_cl}.
\end{proof}

The above theorem allows us to apply general (classical)
pro-representability results to the classical deformation functor
$\Def _{\cl}(E).$ The point is that since the abelian category
$\Qcoh _X$ does not have enough projectives we cannot directly
compare the pseudo-functor $\Def _{\cl}(E)$ to the analogous
deformation pseudo-functor for a DG module over a DG category. But
since $\Qcoh _X$ has enough injectives this can be done for the
pseudo-functor $\coDef _{\cl}(E).$ Namely, we have the following
corollary.

\begin{cor}\label{classical} Let $X$ be a locally noetherian scheme over $k$ and
$E\in D^b_{\Qcoh}(X).$ Choose a bounded below complex $I$ of
injective quasi-coherent sheaves on $X$ which is quasi-isomorphic to
$E.$ Assume that the minimal $A_{\infty}\text{-}$model $A$ of a DG
algebra $\cC=\End(I)$ is admissible (II, Definition 4.1) finite
dimensional Koszul (II, Definition 16.5)
$A_{\infty}\text{-}$algebra. Put $\hat{S}=(B\bar{A})^*,$ where
$B\bar{A}$ is the bar construction of the augmentation
($A_{\infty}\text{-}$)ideal $\bar{A}.$ (Thus $\hat{S}$ is a local
complete DG algebra which is acyclic except in degree zero). Then

a) there exist equivalences of pseudo-functors from $\cart$ to ${\bf
Gpd}$
$$\Def _{\cl}(E)\simeq \coDef _{\cl}(E)\simeq \coDef
_{\cl}(\cC)\simeq \Def _{\cl}(\cC);$$

b) there exists an isomorphism of functors from $\cart $ to $\Set$
$$h_{H^0(\hat{S})}\simeq \pi _0\cdot \Def _{\cl}(E).$$
\end{cor}

\begin{proof} The first and the last equivalences in a) follow from Theorem
\ref{def_cl=codef_cl} and II, Theorem 13.5 respectively. Also b)
follows from II, Theorem 16.7 b) and the middle equivalence in a).
Thus it suffices to prove the equivalence $\coDef _{\cl}(E)\simeq
\coDef _{\cl}(\cC).$ Clearly, $\coDef _{\cl}(E)\simeq \coDef
_{\cl}(I).$

Denote by $\cM$ the abelian category of quasi-coherent sheaves on
$X.$ We may consider $\cM$ as a DG category and then denote by $\cM
^{op}\text{-}mod$ the DG category of DG modules over the opposite DG
category $\cM^{op}.$ Also let $C^{dg}(\cM)$ be the DG category  of
complexes over $\cM.$ Consider the covariant DG functor $h_\bullet
^*:C^{dg}(\cM)\to \cM ^{op}$ defined by
$$h_M^*(N):=\Hom^{\cdot} _{C^{dg}(\cM)}(M,N)^*,$$
where $(\cdot )^*$ denotes the (graded) $k$-dual. Then by Theorem
\ref{compar_codef} this DG functor establishes an equivalence of
pseudo-functors $\coDef _{\cl}(I)\simeq \coDef _{\cl}(h^*_I),$ where
the second pseudo-functor is defined in I, Definitions 10.8, 10.14.
(Notice that the homomorphism of DG algebras $h_I^*:\End (I)\to \End
(h^*_I)$ is a quasi-isomorphism.) Finally, since the DG $\cM
^{op}$-module $h_I^*$ is h-injective (I, Section 3.1) and is bounded
below (I, Definition 11.5) we may apply II, Proposition 9.10 to find
an equivalence of pseudo-functors $\coDef _{\cl}(h_I^*)\simeq \Def
_{\cl}(\cC).$
\end{proof}

\subsection{Explicit description of the equivalence $\Def _{\cl}(E)\simeq \Def _{\cl}(\cC)$}

Assume in the last corollary that the $A_{\infty}\text{-}$algebra
$A$ satisfies the condition (*) in II, Definition 15.3, i.e. the
canonical morphism
$$k\to \bR \Hom _{\bar{A} ^{op}}(\bR \Hom _{\bar{A}}(k,A),A)$$
is a quasi-isomorphism. Then we can make explicit the equivalence
$\Def _{\cl}(E)\simeq \Def _{\cl}(\cC).$

Namely, consider the $A_{\infty}$ $\bar{A}_{\hat{S}}$-module $k.$
It was shown in II, Section 15 that the DG
$(\cC\otimes\hat{S})^{op}$-module
$$\cE=\Hom^{\cdot} _{\bar{A}}(k,\cC)$$
is the universal pro-deformation of the DG $\cC ^{op}$-module $\cC.$
In particular, given a (commutative) local artinian algebra $\cR \in
\cart$ and an object $(T,\tau )\in \Def _{\cR}(\cC)$ there exists a
homomorphism of DG algebras $g:\hat{S}\to \cR$ such that the object
$(\cE \otimes _{g}\cR, \id)$ in $\Def _{\cR}(\cC)$ is isomorphic to
$(T,\tau)$ (it follows from II, Lemma 8.10). (Notice that $\cE$ as a
graded $(\cC \otimes \hat{S})^{op}$-module is isomorphic to $\cC
\otimes \hat{S},$ so actually the graded $\cC^{op}$-module $\cE
\otimes _{g}\cR$ is free of finite rank).

The complex $(\cE \otimes _g\cR)\stackrel{\bL}{\otimes} _{\cC}I$ of
quasi-coherent sheaves on $X_{\cR}$ is an object in $\Def
_{\cR}(I)=\Def _{\cR}(E)$ corresponding to $(T,\tau ).$

\section{Deformation of points objects on a smooth variety and other
examples}

\begin{defi} Let $X$ be a scheme and $E\in D^b(\coh _X).$ We call
$E$ a point object of dimension $d$ if the DG algebra $\bR
\Hom^{\cdot} (E,E)$ is formal, i.e. it is quasi-isomorphic to its
cohomology algebra $\Ext^{\cdot}(E,E),$ and this algebra is
isomorphic to the exterior algebra of dimension  $d.$
\end{defi}

Let $E$ be a point object and put $\cC =\Ext^{\cdot}(E,E).$ By II,
Theorem 15.2 the deformation pseudo-functor $\DEF _-(E)$ is
pro-representable by the DG algebra $\hat{S}=(B\bar{\cC})^*.$ This
DG algebra is quasi-isomorphic to its zero cohomology algebra
$H^0(\hat{S}),$ which is a commutative power series ring. Thus the
formal DG moduli space of point objects is an ordinary (concentrated
in degree zero) commutative regular scheme.

The following proposition justifies our term "point object".

\begin{prop} Let $X$ be a  scheme of finite type over $k$ and let
$p\in X$ be a smooth $k$-point. Then the structure sheaf $\cO _p\in
D^b(\coh _X)$ is a point object of dimension $\dim \cO _{X,p}.$
\end{prop}

\begin{proof} Denote by $j:\Spec \cO _{X,p}\hookrightarrow X$ the
canonical morphism of schemes. It induces an exact functor $j_*:\cO
_{X,p}\text{-mod}\to \Mod _X$ which preserves injective objects
(being the right adjoint to the exact functor $j^*:\Mod _X\to \cO
_{X,p}\text{-mod}$). The functor $j_*$ maps $k$ to $\cO _p$ and
hence induces a quasi-isomorphism of DG algebras
$$j_*:\bR \Hom _{\cO_{X,p}}(k,k) \to \bR \Hom _{D(X)}(\cO _p,\cO
_p). $$ So it suffices to show that the DG algebra $\bR \Hom
_{\cO_{X,p}}(k,k)$ is quasi-isomorphic to the exterior algebra.

Denote the local ring $\cO _{X,p}=A$ and let $m\subset A$ be the
maximal ideal. Consider $A$ as an augmented DG algebra concentrated
in degree zero.

Choose a subspace $V\subset m$ which maps isomorphically to $m/m^2.$
Consider the exterior coalgebra $\wedge ^\bullet V,$ where $\deg
V=-1,$ $\Delta (v)=v\otimes 1+1\otimes v$ for $v\in V$ and $d=0.$
Then the identity map $V\to m$ is an admissible twisting cochain
$\tau \in \Hom^{\cdot} _k(\wedge ^\bullet V,A)$ (II, Definition
2.2). The corresponding DG $A^{op}$-module $\wedge ^\bullet V\otimes
_{\tau}A$ (II, Example 2.6) is just the usual Koszul complex for
$A,$ hence it is quasi-isomorphic to $k.$ Thus
$$\bR \Hom _{A^{op}}(k,k)=\Hom^{\cdot} _{A^{op}}(\wedge ^\bullet V\otimes
_{\tau}A,\wedge ^\bullet V\otimes _{\tau}A).$$

Define a map of complexes
$$\theta :\Hom^{\cdot} _k(\wedge ^\bullet V,k)\to \Hom^{\cdot} _{A^{op}}(\wedge ^\bullet V\otimes
_{\tau}A,\wedge ^\bullet V\otimes _{\tau}A)$$ by the formula $\theta
(f)(a\otimes b)=f(a_{(1)})a_{(2)}\otimes b,$ where $\Delta
(a)=a_{(1)}\otimes a_{(2)}.$ Now exactly as in the proof of II,
Lemma 3.8 one can show that $\theta $ is a homomorphism of DG
algebras, which is a quasi-isomorphism.
\end{proof}

Another example of a point object is provided by a line bundle on an
abelian variety.

On the other hand let $E$ be a line bundle on a (smooth projective)
curve $X$ of genus $g.$ Then the DG algebra $\bR \Hom (E,E)$ is
formal and $\Ext ^0(E,E)=k,$ $\Ext ^1(E,E)=W$ - a vector space of
dimension $g,$ and $\Ext ^i(E,E)=0$ for $i>1.$ By II, Theorem 8.2
the pseudo-functor $\DEF _-(E)$ is pro-representable by the DG
algebra $\hat{S}=(B\overline{\Ext^{\cdot}(E,E)})^*.$ This DG algebra
is concentrated in degree zero and is isomorphic to a noncommutative
power series ring of dimension $g.$

\begin{remark} \rm{The above two examples of line bundles show that the
Picard variety of an abelian variety is (at least locally) the
"true" moduli space of line bundles, whereas the Picard variety of
a curve (of genus $g>1$) is not (!). Indeed, the above argument
shows that in the case of a curve the Picard variety (at least
locally) has a natural noncommutative structure.}
\end{remark}

\part{Noncommutative Grassmanians}

\section{Preliminaries on $\Z$-algebras}
\label{Z-algebras}

In this section we define the notion of a $\Z$-algebra and
associate to it an abelian category which should be thought of as
a category of quasi-coherent sheaves on the corresponding
noncommutative stack. We also define Koszul $\Z$-algebras.

\begin{defi} A $\Z$-algebra $\cA$ over the field $k$ is a
$k$-linear category with the set of objects $\Z.$ For $i,j\in \Z,$
we write $\cA_{ij}$ instead of $\Hom_{\cA}(i,j).$ Sometimes we will
identify a $\Z$-algebra $\cA$ with the corresponding ordinary
non-unital algebra $\Alg_{\cA}.$
\end{defi}

Further, if $\cA$ is a $\Z$-algebra, then we define the abelian
category $\text{Mod-}\cA$ as the category
$\Fun(\cA^{op},k\text{-Vect})$ of contravariant functors from
$\cA$ to $k$-vector spaces. Equivalently $\text{Mod-}\cA$ is the
full subcategory of $\text{Mod-}\Alg_{\cA}$ which consists of
right $\Alg_{\cA}$-modules $M$ such that
$$M=\bigoplus\limits_{i\in\Z}M\cdot\one_i$$ (quasi-unital modules). We
call the objects of $\text{Mod-}\cA$ $\cA^{op}$-modules. For each
$i\in\Z$ put
$$P_i:=\Hom(-,i)=\one_i Alg_{\cA}\in \text{Mod-}\cA.$$ By Ioneda Lemma, for each $M\in \text{Mod-}\cA$ we have
$$\Hom_{\cA^{op}}(P_i,M)=M(i),$$ hence $P_i$ are projectives. Clearly, each
$M\in \text{Mod-}\cA$ can be covered by a direct sum of $P_i$'s,
hence the abelian category $\text{Mod-}\cA$ has enough
projectives.

\begin{defi} Let $M\in \text{Mod-}\cA$ be an $\cA^{op}$-module. An element $x\in M(i)$ is
called torsion if we have $x\cA_{ji}=0$ for $j<<i.$ Torsion elements
form a submodule of $M$ which we denote by $\tau(M).$ An
$\cA^{op}$-module $M$ is called torsion if we have $M=\tau(M).$ We
denote by $\Tors(\cA)$ the full subcategory of $\text{Mod-}\cA$
which consists of torsion $\cA^{op}$-modules.

The category $\QMod(\cA)$ is defined as the quotient category
$\text{Mod-}\cA/\Tors(\cA).$ We denote by $\pi:\text{Mod-}\cA\to
\QMod(\cA)$ the projection functor.
\end{defi}

If $M,$ $N$ are $\cA^{op}$-modules then
$$\Hom_{\QMod(\cA)}(\pi(M),\pi(N))=\lim_{\rightarrow}\Hom_{\cA^{op}}(M',N/\tau(N)),$$
where $M'$ runs over the quasi-directed category of submodules
$M'\subset M$ such that $M/M'$ is torsion.

The category $\QMod(\cA)$ should be thought of as the category
$\QCoh(\Proj(\cA))$ of quasi-coherent sheaves on the noncommutative
projective stack $\Proj(\cA).$ Furthermore, the object $\pi(P_0)\in
\QMod(\cA)$ should be thought of as the structure sheaf
$\cO_{\Proj(A)}.$

\begin{remark}\label{graded-algebras} \rm{Let $A=\bigoplus\limits_{i\in\Z}A^i$ be a
$\Z$-graded (unital) algebra. Then one can associate to it a
$\Z$-algebra $\cA$ with $\cA_{ij}=A^{j-i}$ so that the composition
in $\cA$ comes from the multiplication in $A.$ Recall that in
\cite{AZ} there defined the category $\QGr(A)$ as the quotient
category $\Gr A/\Tors$ of the category $\Gr A$ of graded $A$-modules
by the subcategory $\Tors$ of torsion modules. It is clear that the
categories $\Gr A,$ $\Tors,$ $\QGr A$ are equivalent  to
$\text{Mod-}\cA,$ $\Tors(\cA),$ $\QMod(\cA)$ respectively.

Notice that it can happen that the graded algebras $A_1$ and $A_2$
are not isomorphic but the associated $\Z$-algebras are
equivalent. Thus, it is more reasonable to consider $\Z$-algebras
but not graded algebras.}
\end{remark}

The projection functor $\pi:\text{Mod-}\cA\to \QMod(\cA)$ admits a
right adjoint functor $\omega:\QMod(\cA)\to \text{Mod-}\cA$
defined by the formula $$\omega(X)(i)=\Hom(\pi(P_i),X).$$ The
adjunction morphism $\pi\omega\to \id$ is an isomorphism.

\begin{defi} A $\Z$-algebra $\cA$ is called

a) positively (resp. negatively) oriented if $\cA_{ij}=0$ for
$i>j$ (resp. for $i<j$);

b) connected if $\cA_{ii}\cong k$ for each $i\in\Z$;

c) locally finite if $\dim\cA_{ij}<\infty$ for any $i,j\in\Z.$
\end{defi}

Let $\cA$ be a positively oriented $\Z$-algebra. We denote by
$\cA_{\leq i}$ the full subcategory of $\cA$ such that $Ob(A_{\leq
i})=\{j: j\leq i\}.$ Clearly, we also have the categories
$\text{Mod-}\cA_{\leq i}$ and $\Tors(\cA_{\leq i}).$ It is easy to
see that the quotient category $\text{Mod-}\cA_{\leq
i}/\Tors(\cA_{\leq i})$ is equivalent to $\QMod(\cA).$ We denote by
$\pi_{\leq i}:\text{Mod-}\cA_{\leq i}\to \QMod(\cA)$ and
$\omega_{\leq i}:\QMod(\cA)\to \text{Mod-}\cA_{\leq i}$ the
projection functor and its right adjoint respectively.

If $\cA$ is a positively oriented $\Z$-algebra then we put
$$T_{ij}=P_j/(P_j)_{\leq i},$$ where $$(P_j)_{\leq i}=\bigoplus\limits_{k\leq
i}A_{kj}.$$ Clearly, the $\cA^{op}$-modules $T_{ij}$ are torsion.

If $\cA$ is a positively or negatively oriented connected
$\Z$-algebra then we denote by $S_n$ the simple $\cA^{op}$-modules
defined by the formula
$$S_n(i)=\begin{cases}k & \text{for }i=n,\\
0 & \text{otherwise.}\end{cases}$$ Notice that if $\cA$ is
positively oriented then $S_n=T_{n,n}.$

\begin{defi} A connected positively (resp. negatively) oriented $\Z$-algebra
is called quadratic if the algebra $Alg_{\cA}$ is generated by the
subspaces $\cA_0=\bigoplus\limits_{i\in \Z} \cA_{ii}$ and
$\cA_1=\bigoplus\limits_{i\in \Z} \cA_{i,i+1}$ (resp.
$\cA_{-1}=\bigoplus\limits_{i\in \Z} \cA_{i+1,i}$) and is
determined by the quadratic relations $I_{i,i+2}\in
\cA_{i+1,i+2}\otimes \cA_{i,i+1}$ (resp. $I_{i+2,i}\in
\cA_{i+1,i}\otimes \cA_{i+2,i+1}$).
\end{defi}

For a locally finite positively (resp. negatively) oriented
quadratic $\Z$-algebra $\cA$ one can define the dual quadratic
$\Z$-algebra $\cA^{!}$ with the opposite orientation. It is
defined by the dual generators $\cA^!_{i+1,i}=\cA_{i,i+1}^*$
(resp. $\cA^!_{i,i+1}=\cA_{i+1,i}^*$) and the dual quadratic
relations $S(I_{i,i+2}^{\perp})\subset \cA_{i,i+1}^*\otimes
\cA_{i+1,i+2}^*$ (resp. $S(I_{i+2,i}^{\perp})\subset
\cA_{i+2,i+1}^*\otimes \cA_{i+1,i}^*$), where
$I_{i,i+2}^{\perp}\subset \cA_{i+1,i+2}^*\otimes \cA_{i,i+1}^*$
(resp. $I_{i+2,i}^{\perp}\subset \cA_{i+1,i}^*\otimes
\cA_{i+2,i+1}^*$) is the dual subspace and $S$ is the
transposition of factors.

Further, one can define a Koszul complex
$$K:=\cA^{!^*}\otimes_{\cA_0}\cA=\bigoplus \cA^{!^*}_{kj}\otimes_k \cA_{ij}.$$
Here $\cA^{!^*}=\bigoplus\limits_{i,j} \cA^{!^*}_{ij}$ is a bounded
dual of $\cA^{!}.$ It is an $\cA^{!}$-bimodule.

The differential $d:K\to K$ is defined as follows. Suppose that
$\cA$ is positively oriented. We have the natural maps
$$\cA_{j,j+1}\otimes \cA_{ij}\to \cA_{i,j+1},\quad\text{and}\quad
\cA_{j,j+1}^*\otimes \cA_{kj}^{!^*}\to \cA_{k,j+1}^{!^*}.$$ In
particular, we have the maps $$\psi_{ijk}:\cA_{j,j+1}^*\otimes
\cA_{j,j+1}\to \Hom_k (\cA^{!^*}_{kj}\otimes_k
\cA_{ij},\cA^{!^*}_{k,j+1}\otimes_k \cA_{i,j+1}).$$ The non-zero
components of $d$ are the maps
$d_{ijk}=\psi_{ijk}(\one_{A_{jj+1}}).$ Note that $d$ is
$\cA_0$-linear and $\cA^!\otimes \cA^{op}$-linear. Thus, $K_n=\one_n
K$ and $K_n^m=K_n\one_m$ are $d$-invariant. The complex $K_n$ is of
the form $$\dots\to \cA^{!^*}_{n,n-2}\otimes P_{n-2}\to
\cA^{!^*}_{n,n-1}\otimes P_{n-1}\to P_n\to 0,$$ and the complex
$K_n^m$ is of the form $$\dots\to \cA^{!^*}_{n,n-2}\otimes
\cA_{m,n-2}\to \cA^{!^*}_{n,n-1}\otimes \cA_{m,n-1}\to \cA_{m,n}\to
0.$$ In particular, $K_n^n\cong k.$

Analogously for negatively oriented $\Z$-algebras.

For the rest of this section we assume that $\cA$ is positively
oriented.

\begin{defi} A quadratic locally finite $\Z$-algebra is called Koszul if the complex $K_n^m$
is acyclic for $n\ne m,$ or, equivalently, $K_n$ is a resolution of
$S_n.$\end{defi}

We refer to [BP] for the definition of co-Koszul and Gorenstein
$\Z$-algebras. We will not need these definitions but we will need
the following proposition:

\begin{prop}(\cite{BP})\label{Frobenious} Let $\cA$ be a Koszul (positively oriented) $\Z$-algebra of finite
homoogical dimension $n.$ Then the following conditions are
equivalent:

(i) $\cA$ is co-Koszul;

(ii) $\cA$ is Gorenstein;

(iii) $\cA^!$ is Frobenious, i.e. $\cA_{i+n,i}^!=k$ for all $i,$ and
the multiplication $\cA_{ji}\otimes \cA_{i+n,j}\to k$ is a perfect
pairing for all $i$ and $j.$
\end{prop}

Now we define the notion of a coherent $\Z$-algebra and the category
$\qmod(\cA)$ for a coherent $\Z$-algebra $\cA.$

\begin{defi} Let $\cA$ be a $\Z$-algebra. An $\cA^{op}$-module $M$
is called finitely generated if there exists a surjective morphism
(in $\text{Mod-}\cA$)
$$\bigoplus\limits_{j=1}^m P_{i_j}\to M,$$ where $i_1,\dots,i_m\in
\Z.$ Further, a finitely generated $\cA^{op}$-module $M$ is called
coherent if for each (not necessarily surjective) morphism (in
$\text{Mod-}\cA$) $$\phi:\bigoplus\limits_{j=1}^m P_{i_j}\to M$$ the
$\cA^{op}$-module $\ker(\phi)$ is finitely generated. A $\Z$-algebra
$\cA$ is called coherent if for each $i\in\Z$ the module $P_i$ is
coherent.

If $\cA$ is coherent then we denote by $\qmod(\cA)$ the full
(abelian) subcategory of $\QMod(\cA)$ which consists of the images
of coherent $\cA^{op}$-modules.
\end{defi}

The category $\qmod(\cA)$ should be thought of as the category
$\Coh(\Proj(\cA)).$ By definition, we have that the structure sheaf
$\cO_{\Proj(\cA)}$ is coherent.

\section{The definition of noncommutative Grassmanians}
\label{Grassmans}

Let $V$ be a finite-dimensional $k$-vector space of dimension $n>0.$
Let $m$ be an integer such that $1\leq m\leq n-1.$ We define the
noncommutative Grassmanians by the formula
$$\NGr(m,V):=\Proj(\cA^{m,V}),$$ where $\cA^{m,V}$ is the following quadratic $\Z$-algebra:
$$\cA^{m,V}_{i,i+1}=\begin{cases}V^* & \text{for }(n-m+1)\nmid i,\\
\Lambda^{n-m}V & \text{otherwise,}\end{cases}$$ and the quadratic
relations are defined by the natural exact sequences
$$\begin{cases}\Lambda^2 V^*\to\cA^{m,V}_{i+1,i+2)}\otimes\cA^{m,V}_{i,i+1}\to\cA^{m,V}_{i,i+2}\to 0 & \text{for }
(n-m+1)\nmid i,i+1,\\
\Lambda^{n-m-1}V\to\cA^{m,V}_{i+1,i+2}\otimes\cA^{m,V}_{i,i+1}\to\cA^{m,V}_{i,i+2}\to
0 & \text{otherwise.}
\end{cases}
$$

Notice that if we fix a volume form $\omega\in\Lambda^n V,$ then the
$\cA^{1,V}$ is naturally equivalent to the $\Z$-algebra associated
to the symmetric algebra $\bigoplus\limits_{l\geq 0}S^l V^*,$ where
$\deg(V^*)=1.$ Hence, the stack $\NGr(1,V)$ is isomorphic to the
commutative projective space $\PP(V).$

We claim that $\NGr(m,V)$ is a true noncommutative moduli space of
structure sheaves $\cO_{\PP(W)}\in D^b_{coh}(\PP(V)),$ where
$W\subset V$ are vector subspaces of dimension $\dim W=m.$ Namely,
it satisfies the following properties:

1) There is a natural fully faithful functor $\Phi$ from the
category of perfect objects
$\Perf(\NGr(m,V))$ (Definition
\ref{perfect}) to $D^b_{coh}(\PP(V)).$ Its image is the double
orthogonal to the family of objects $\cO_{\PP(W)},$ i.e. the full
subcategory generated by objects
$\cO_{\PP(V)}(m-n),\dots,\cO_{\PP(V)}(-1),\cO_{\PP(V)}.$ This is
Corollary \ref{T^m,V} below;

2) There is a $k$-point $x_W\in \NGr(m,V)(k)=X_{\cA^{m,V}}(k)$
(see Section \ref{k-points} below) for each subspace $W\subset V$ of
dimension $\dim W=m.$ Further, $(x_W)_*(k)$ lies in
$\Perf(\NGr(m,V))$ and $\Phi(x_*(k))\cong \cO_{\PP(W)}.$ This is a
part of Theorem \ref{all_k-points} below.

3) The completion of the local ring of the $k$-point $x_W$ (see
Section \ref{defi_local}) is isomorphic to $H^0(\hat{S}),$ where
$\hat{S}$ is dual to the bar construction of the minimal
$A_{\infty}$-structure on $\Ext^{\cdot}(\cO_{\PP(W)},\cO_{\PP(W)})$
(Theorem \ref{local_rings}).  It can be shown that the DG algebra
$\bR\Hom(\cO_{\PP(W)},\cO_{\PP(W)})$ is formal and the graded
algebra $\Ext^{\cdot}(\cO_{\PP(W)},\cO_{\PP(W)})$ is quadratic
Koszul, and hence the projection $\hat{S}\to H^0(\hat{S})$ is a
quasi-isomorphism. Hence, the moduli space is not a DG space but
just noncommutative space.

Furthermore, we do not have a moduli functor of our family of
objects $\cO_{\PP(W)},$ which should be defined on the category of
noncommutative affine schemes. However, the properties 1), 2), and 3)
suggest that $\NGr(m,V)$ is a true moduli space of this family of
objects, in our context of deformations of objects in derived
categories.

It is remarkable that there is a natural morphism from the
commutative Grassmanian $\Gr(m,V)$ to noncommutative one
$\NGr(m,V).$ Moreover, the functor $\Phi:\Perf(\NGr(m,V))\to
D^b_{coh}(\PP(V))$ above coincides with $\bL f_{1,m,V}^*,$ where
$f_{1,m,V}:\PP(V)\to \NGr(m,V)$ is a natural morphism. Both these
statements are parts of Proposition \ref{Gr_to_NGr} below.

\section{The derived categories of noncommutative Grassmanians}
\label{derived_on_NGr}

Before we formulate and prove results on the derived categories of
quasi-coherent sheaves on noncommutative Grassmanians $\NGr(m,V)$
we need to remind some notions and results from \cite{BP}.

Let $D$ be a $k$-linear enhanced triangulated category.

\begin{defi} An object $E\in Ob (D)$ is called exceptional if
$\Hom^i(E,E)=0$ for $i\ne 0,$ and $\Hom^0(E,E)=k.$
\end{defi}

\begin{defi} A collection $(E_1,\dots,E_m)$ of exceptional objects in $D$ is called exceptional
if $\Hom^*(E_i,E_j)=0$ for $i>j.$\end{defi}

\begin{defi} A full exceptional collection of objects in the category $D$ is a collection which
generates $D$ as triangulated category.\end{defi}

\begin{defi} An exceptional collection $(E_1,\dots,E_n)$ is called strong exceptional if
it satisfies the additional assumption $\Hom^i(E_k,E_l)=0$ for $i\ne
0$ and all $k$ and $l.$
\end{defi}

Let $(E,F)$ be an exceptional pair. Define the objects $L_{E}F$
and $R_{F}E$ by the exact triangles

$$L_E F\to \Hom^{\cdot}(E,F)\otimes E\to F,\quad E\to \Hom^{\cdot}(E,F)^{\vee}\otimes F\to R_F E.$$

Let $\sigma=(E_1,\dots,E_n)$ be an exceptional collection. If
$1\leq i\leq n-1$ (resp. $2\leq i\leq n$), then the right (resp.
left) mutation of the object $E_i$ in this collection is the
object $R^1 E_i=R_{E_{i+1}}E_i$ (resp. $L^1 E_i=L_{E_{i-1}}E_i$);
the corresponding mutated collection
$$R^1_{E_i}\sigma:=(E_1,\dots,E_{i-1},E_{i+1},R_{E_{i+1}}E_i,E_{i+2},\dots,E_n)$$
(and the analogous collection $L^1_{E_i}\sigma$) is exceptional.
The multiple mutations of the objects and of the collection are
defined inductively:

$$R^k E_i=R_{E_{i+k}}R^{k-1}E_i,\quad R^k_{E_i}\sigma= R^1_{R^{k-1}E_i}(R^{k-1}_{E_i}\sigma),\quad k\leq n-i$$
(and in the same way for left mutations).

\begin{defi} A helix of the period $n$ is an infinite sequence $\{E_i\}_{i\in \Z}$
such that for each $i\in \Z$ the collection $(E_i,\dots,E_{i+n})$ is
exceptional, and moreover $R^{n-1}E_i=E_{i+n}.$\end{defi}

If $\sigma=(E_1,\dots,E_n)$ is an exceptional collection then it
naturally extends to a helix by the conditions
$$E_{i+n}=R^{n-1}E_i,\quad i\geq 1,\quad\quad E_{i-n}=L^{n-1}E_i,\quad i\leq n.$$
In this case the helix is said to be generated by the collection
$\sigma.$

If the helix is generated by the full exceptional collection then it
satisfies the property of the partial periodicity: $\Phi (E_i)\cong
E_{i-n},$ where $\Phi=F[1-n]$ is the composition of the Serre
functor $F$ and the multiple shift $[1-n].$

\begin{defi}(\cite{BP}) A helix $\mu=\{E_i\}$ is called geometric if for each pair $(i,j)\in\Z^2$
such that $i\leq j$ one has
$$\Hom^k(E_i,E_j)=0\text{ for }k\ne 0.$$\end{defi}

\begin{defi}(\cite{BP}) An exceptional collection is called geometric if it generates a geometric helix.\end{defi}

\begin{prop}(\cite{BP})\label{sub-collection} Each sub-collection of a geometric exceptional collection is again geometric.\end{prop}

\begin{prop}(\cite{BP})\label{Fano} A full exceptional collection of the length $m$ of coherent sheaves
on a smooth projective variety $X$ of dimension $n$ is geometric iff
$m=n+1.$\end{prop}

\begin{defi}(\cite{BP}) The endomorphism $\Z$-algebra $\cA=\End(S)$ of a helix $S=\{E_i\}$ is defined by the formula
$$\cA_{ij}=\Hom(E_i,E_j)$$ with natural composition.\end{defi}

\begin{theo}(\cite{BP})\label{End(S)_geom} If $S$ is a geometric helix generated by an exceptional collection
of length $n$ then the endomorphism $\Z$-algebra $\cA$ of $S$ is
Koszul, co-Koszul and of finite global homological dimension $n.$
\end{theo}

\begin{defi}(\cite{BP}) Koszul co-Koszul $\Z$-algebra of finite homological dimension $n$ is
called a geometric $\Z$-algebra of the period $n.$\end{defi}

Let $\cA$ be a geometric $\Z$-algebra of the period $n.$ Let
$K\subset D(\text{Mod-}\cA)$ be the full triangulated subcategory
generated by the modules $P_i.$ Note that $S_i\in K.$ Let $F\in K$
be the full triangulated subcategory generated by the modules $S_i.$

\begin{theo}(\cite{BP})\label{K/F} Let $\cA$ be a geometric $\Z$-algebra of the period $n.$ Then
$F$ is a thick subcategory in $K$; the images of modules $P_i$ in
$K/F$ form a geometric helix $S$ of the period $n,$ and moreover the
$\Z$-algebra of $S$ is equivalent to $\cA.$
\end{theo}

Now we prove the main theorem of this section. It is closely
related to the previous one but unfortunately cannot be deduced
from it.

\begin{theo}\label{D(Geom)} Let $\cA$ be a geometric helix of the period $n.$ Put $B=\cA_{[1,n]}=
\bigoplus\limits_{1\leq i,j\leq n} \cA_{ij}.$ Then there is an
equivalence of categories $D^*(\QMod(\cA))\cong D^*(\text{Mod-}B).$
\end{theo}
\begin{proof} The proof is in two main steps. First we prove that
the category $D^*(\QMod(\cA))$ is naturally equivalent to the
quotient category of $D^*(\text{Mod-}\cA_{\leq n})$ by the full
thick triangulated subcategory $D^*_{\Tors}(\cA_{\leq n})$ which
consists of complexes with torsion cohomology. Then we construct
mutually inverse exact equivalences between the categories
$D^*(\text{Mod-}\cA_{\leq n})/D^*_{\Tors}(\cA_{\leq n})$ and
$D^*(\text{Mod-}\cB)$ given by DG bimodules.

\begin{lemma}\label{quotient} The categories $D^*(\QMod(\cA))$ and $D^*(\text{Mod-}\cA_{\leq n})/D^*_{\Tors}(\cA_{\leq n})$
are naturally equivalent.\end{lemma}

\begin{proof} First recall the functor $\omega_{\leq n}:\QMod(\cA)\to
\text{Mod-}\cA_{\leq n}.$ It induces a fully faithful functor
$$K^*(\omega_{\leq n}):K^*(\QMod(\cA_{\leq n}))\to K^*(\text{Mod-}\cA_{\leq n})$$
between homotopy categories, which is right adjoint to the functor
$K^*(\pi_{\leq n}).$ It follows that $K^*(\QMod(\cA_{\leq n}))$ is
equivalent to the quotient category $K^*(\text{Mod-}\cA_{\leq
n})/K^*(\pi_{\leq n})^{-1}(0).$ Let $K^*_{\Tors}(\cA_{\leq
n})\subset K^*(\text{Mod-}\cA_{\leq n})$ be the full subcategory
that consists of all complexes with torsion cohomology. It is easy to see that acyclic
complexes in the category $K^*(\QMod(\cA_{\leq n}))$ correspond to
the classes of complexes with torsion cohomology in
$K^*(\text{Mod-}\cA_{\leq n})/K^*(\pi_{\leq n})^{-1}(0).$ Thus,
$D^*(\QMod(\cA_{\leq n}))$ is equivalent to the quotient of
$K^*(\text{Mod-}\cA_{\leq n})/K^*(\pi_{\leq n})^{-1}(0)$ by
$K^*_{\Tors}(\cA_{\leq n})/K^*(\pi_{\leq n})^{-1}(0).$ This quotient
is further equivalent to $D^*(\text{Mod-}\cA_{\leq
n})/D^*_{\Tors}(\cA_{\leq n}).$ Lemma is proved.
\end{proof}

Denote by $(Q_1,\dots,Q_n)$ the exceptional collection of
indecomposable projective $B^{op}$-modules. By Theorem \ref{K/F} the
helix $\{Q_i\}_{i\in \Z}$ generated by $(Q_1,\dots,Q_n)$ is
geometric. It follows from its partial periodicity that for $i\leq
0$ we have $Q_i\cong H^{n-1}(Q_i)[1-n].$ Thus, we may and will
assume that $Q_i$ is concentrated in degree $n-1$ for $i\leq 0.$ Put
$M_0=\bigoplus\limits_{i\leq 0} Q_i[n-1].$ Since
$\Hom_B(Q_i,Q_j)=\cA_{ij}=\cA_{(i+n)(j+n)},$ $M_0$ is naturally an
$\cA_{\leq n}\otimes B^{op}$-module. Further, the functor
$\Phi^{-1}=F^{-1}[1-n]$ can be given by the object $B^{!}[1-n]\in
D(B\otimes B^{op}),$ where $B^{!}=\bR \Hom_{B^{op}\otimes
B}(B,B\otimes B).$ Since $Q_i\otimes^{\bL}_{B} B^{!}[1-n]$ are pure
modules for $i=1,\dots,n,$ $B^{!}[1-n]$ is a pure bimodule. We
define the object $M_1\in D^b(\cA_{\leq n}\otimes B^{op})$ by the
formula
$$M_1:=M_0\stackrel{\bL}{\otimes}_{B} B^{!}[2-2n].$$ We have that
$$P_i\stackrel{\bL}{\otimes}_{\cA_{\leq n}} M_1\cong Q_i$$
for $i\leq n.$ Since $\cA_{\leq n}$ has finite left homological
dimension, we have a well defined  functor
$$-\stackrel{\bL}{\otimes}_{\cA_{\leq n}} M_1:D^*(\text{Mod-}\cA_{\leq n})\to
D^*(\text{Mod-}B).$$

\begin{lemma}\label{torsion_to_0} For each $K^{\cdot}\in D_{\Tors}(\cA_{\leq n}^{op})$ we have
$$K^{\cdot}\stackrel{\bL}{\otimes}_{\cA_{\leq n}}M_1=0.$$\end{lemma}
\begin{proof} Clearly, it suffices to prove the Lemma for
$M_0$ instead of $M_1.$ We have that the complex $K_l$ from section
\ref{Z-algebras} is a projective resolution of $S_l$ for $l\leq n.$
Further, $S_l\stackrel{\bL}{\otimes}_{\cA_{\leq n}}M_0\cong
K_l\otimes_{\cA_{\leq n}}M_0$ and the last complex is up to shift of
the following form

$$\dots0\to A_{l-n,l-2n}^{!^*}\otimes
Q_{l-2n}[n-1]\to\dots\to A_{l-n,l-n-1}^{!^*}\otimes
Q_{l-n-1}[n-1]\to Q_{l-n}[n-1]\to 0\to \dots$$

This complex is acyclic since it corresponds to the image of
$K_{l-n}$ in $K/F$ under the equivalence of Theorem \ref{K/F}, and
the image of $K_{l-n}$ in $K/F$ is zero.

Further, the torsion modules $T_{km},$ $k\leq m\leq n$ have finite
filtrations with subquotients being direct sums of $S_l.$ Thus, we
have $T_{km}{\otimes}^{\bL}_{\cA_{\leq n}}M_0=0.$

Since each torsion $\cA_{\leq n}^{op}$-modules has a left resolution
by the direct sums of $T_{km},$ it follows that the statement of the
Lemma holds if $K^{\cdot}$ is a pure torsion $\cA_{\leq
n}^{op}$-module. Finally, since $M_0$ is quasi-isomorphic isomorphic
to a finite complex of bimodules which are projective as left
$\cA$-modules, the statement of the Lemma holds for each
$K^{\cdot}\in D_{\Tors}(\cA_{\leq n}^{op}).$
\end{proof}

By the previous Lemma, the formula
$-\stackrel{\bL}{\otimes}_{\cA_{\leq n}} M_1$ defines a functor
$$\Phi:D^*(\text{Mod-}\cA_{\leq n})/D^*_{\Tors}(\cA_{\leq n})\to D^*(\text{Mod-}B).$$

Further, $M_2:=\bigoplus\limits_{1\leq i\leq n} P_i$ is naturally
an $\cA_{\leq n}^{op}\otimes B$-module. Consider the functor
$$\Psi:D^*(\text{Mod-}B)\to D^*(\text{Mod-}\cA_{\leq
n})/D^*_{\Tors}(\cA_{\leq n})$$ defined by the formula
$$\Psi(-)=\pi_{tors}(-\stackrel{\bL}{\otimes}_B M_2),$$
where $\pi_{tors}:D^*(\text{Mod-}\cA_{\leq n})\to
D^*(\text{Mod-}\cA_{\leq n})/D^*_{\Tors}(\cA_{\leq n})$ is the
projection.

\begin{lemma}\label{inverse} The functors $\Phi$ and $\Psi$ are mutually inverse equivalences.\end{lemma}
\begin{proof} First, the isomorphism
$$M_2\stackrel{\bL}{\otimes}_{\cA_{\leq n}} M_1\to B$$ in $D(B\otimes B^{op})$
induces the isomorphism of functors $\Phi\circ\Psi\cong Id.$

Further, we claim that $H^0(M_1\otimes_B^{\bL} M_2)\cong \cA_{\leq
n},$ $H^{n-1}(M_1\otimes_B^{\bL} M_2)$ is torsion as
$\cA^{op}$-module and $H^i(M_1\otimes_B^{\bL} M_2)=0$ for $i\ne
0,n-1.$ Indeed, since $K_l$ is a resolution of $S_l$ it follows from
Lemma \ref{torsion_to_0} by decreasing induction on $l\leq n$ that
$$H^i(P_l\stackrel{\bL}{\otimes}_{\cA_{\leq n}} M_1\stackrel{\bL}{\otimes}_B M_2)=
\begin{cases}P_l & \text{for }i=0;\\
\text{is torsion} & \text{for }i\geq n-1\\
0 & \text{otherwise.}\end{cases}$$ Thus, it remains to note that
$H^{k}(M_1)=0$ for $k\geq n$ and $M_2$ is pure.

Finally, we have the natural morphism $\cA_{\leq n}\to
M_1\otimes_B^{\bL} M_2$ in $D(A_{\leq n}\otimes A_{\leq n}^{op})$
and for each $K^{\cdot}\in D(\text{Mod-}\cA_{\leq n})$ we have that
$$Cone(K^{\cdot}\to K^{\cdot}\stackrel{\bL}{\otimes}_{\cA_{\leq n}}
M_1\stackrel{\bL}{\otimes}_B M_2)\cong
K^{\cdot}\stackrel{\bL}{\otimes}_{\cA_{\leq n}}Cone(\cA_{\leq n }\to
M_1\stackrel{\bL}{\otimes}_B M_2)\in D_{\Tors}(\cA_{\leq n}^{op}).$$
Thus, $\Psi\circ \Phi\cong Id.$ Lemma is proved.
\end{proof}
Theorem follows from Lemmas \ref{inverse} and \ref{quotient}.
\end{proof}

Now we apply the above theorem to noncommutative Grassmanians
introduced in section \ref{Grassmans}. By Propositions
\ref{sub-collection} and \ref{Fano} we have that the exceptional
collection
$$\sigma=(\cO_{\PP(V)}(m-n),\dots,\cO_{\PP(V)}(-1),\cO_{\PP(V)})$$
of coherent sheaves on $\PP(V)$ is geometric. Let $S=\{E_i\}$ be the
helix generated by $\sigma,$ so that $E_i=\cO_{\PP(V)}(i)$ for
$i=m-n,\dots,-1,0.$

\begin{prop}\label{A^m,V=End(S)} The endomorphism
$\Z$-algebra $\cA$ of the helix $S$ is equivalent to
$\cA^{m,V}.$\end{prop}

\begin{proof} Note that both $\cA$ and $\cA^{m,V}$ are quadratic
and $(n-m+1)$-periodic. It remains to show that the space
$\cA_{i,i+1}$ is isomorphic to $\cA^{m,V}_{i,i+1}$ for
$i=m-n,\dots,-1,0,$ and the quadratic relations $I_{i,i+2}\in
\cA_{i+1,i+2}\otimes \cA_{i,i+1}$ coincide with that of $\cA^{m,V}$
for $i=m-n,\dots,0.$

All of this is clear for $i=m-n,\dots,-1.$ Further, the object
$E_1$ is isomorphic to the complex
$$\dots\to0\to\cO_{\PP(V)}(m-n)\to V\otimes \cO_{\PP(V)}(m-n+1)\to\dots\to\Lambda^{n-m}V\otimes \cO_{\PP(V)}\to 0\to\dots,$$
where the last non-zero term is in degree zero. It follows that
$$\cA_{0,1}=\Hom(E_0,E_1)=\Lambda^{n-m} V=\cA^{m,V}_{0,1}.$$
Furthermore, the quadratic relation $I_{-1,1}\subset
\cA_{0,1}\otimes \cA_{-1,0}$ coincides with the subspace
$\Lambda^{n-m-1} V\subset \Lambda^{n-m}V\otimes V^*,$ as in
$\cA^{m,V}.$ Finally, $E_2$ is the convolution of the complex
$$\dots\to 0\to E_{m-n+1}\to V\otimes E_{m-n+2}\to\dots\to \Lambda^{n-m-1}V\otimes E_0\to V^*\otimes E_1\to 0\to \dots,$$
where the last non-zero term is in degree zero. It follows that the
quadratic relation $I_{0,2}\subset \cA_{1,2}\otimes \cA_{0,1}$
coincides with the subspace $\Lambda^{n-m-1} V\subset V^*\otimes
\Lambda^{n-m} V,$ as in $\cA^{m,V}.$
\end{proof}

Let $B^{m,V}$ be the endomorphism algebra of $\sigma.$ As a
corollary of the above results we obtain the following

\begin{theo}\label{D(NGr)} The derived category $D^*(\QMod(\cA^{m,V}))$ is equivalent to
the derived category $D^*(B^{m,V}).$ The objects $\pi(P_i)$ in
$D^*(\QMod(\cA^{m,V}))$ form a geometric helix of the period $\dim
V-m+1.$
\end{theo}
\begin{proof} Indeed, by Proposition \ref{A^m,V=End(S)} and Theorem \ref{End(S)_geom} the $\Z$-algebra $\cA^{m,V}$ is geometric.
Thus, the first statement follows from Theorem \ref{D(Geom)}.
After that, the second statement follows from Theorem
\ref{K/F}.\end{proof}

Now we introduce the perfect derived category.

\begin{defi}\label{perfect} Let $\cA$ be a positively oriented $\Z$-algebra. The category of perfect objects
$\Perf(\QMod(\cA))$ is the minimal full thick triangulated
subcategory of $D(\QMod(\cA))$ which contains the objects
$\pi(P_i).$ We will call the objects of $\Perf(\QMod(\cA))$
perfect complexes.
\end{defi}

We will also write below $\Perf(\NGr(m,V))$ instead of
$\Perf(\QMod(\cA^{m,V})).$

\begin{prop}\label{D_perf} Let $\cA$ be a geometric $\Z$-algebra of the period
$n,$ and $B=\bigoplus_{1-n\leq i,j\leq 0}\cA_{ij}.$ Then the
category $\Perf(\QMod(\cA))$ is equivalent to
$D^b(\text{mod}_{finite}\text{-}B).$\end{prop}
\begin{proof} By Theorem \ref{D(Geom)}
$(\pi(P_{1-n}),\dots,\pi(P_{-1}),\pi(P_0))$ is full strong
exceptional collection in $\Perf(\QMod(\cA)).$ Further, the
category $\Perf(\QMod(\cA))$ is enhanced and
$\End(\bigoplus\limits_{1-n\leq i\leq 0}\pi(P_i))=B.$ Hence
$\Perf(\QMod(\cA))$ is equivalent to
$D^b(\text{mod}_{finite}\text{-}B).$
\end{proof}

\begin{cor}\label{T^m,V} The category $\Perf(\QMod(\cA^{m,V}))$ is
equivalent to the full triangulated subcategory $\T^{m,V}\in
D^{b}_{coh}(\PP(V))$ generated by the exceptional collection
$$(\cO_{\PP(V)}(m-n),\dots,\cO_{\PP(V)}(-1),\cO_{\PP(V)}).$$ Under
this equivalence the exceptional collection
$(\pi(P_{m-n}),\dots,\pi(P_0))$ corresponds to the exceptional
collection
$(\cO_{\PP(V)}(m-n),\dots,\cO_{\PP(V)}(-1),\cO_{\PP(V)}).$
\end{cor}

\begin{remark} \rm{Notice that by \cite{Mi} the $\Z$-algebra $\cA^{\dim
V-1,V}$ is coherent. Further, the $\Z$-algebra $\cA^{1,V}$ is
Noetherian and hence is coherent. It should be plausible that all
the $\Z$-algebras $\cA^{m,V}$ (and, more generally, all geometric
$\Z$-algebras) are coherent, but it is not clear how to prove this
statement. However, if $\cA$ and $B$ are as in Theorem
\ref{D(Geom)}, and $\cA$ is coherent, then the subcategory
$\Perf(\QMod(\cA))\subset D(\QMod(\cA))$ coincides with the
subcategory $D^b_{\qmod}(\QMod(\cA))$ which consists of complexes
with cohomology lying in $\qmod.$ This category is further
equivalent to $D^b(\qmod(\cA)).$ Therefore, in this case we also
have an  equivalence
$$D^b(\qmod(\cA))\cong D^b(\text{mod}_{finite}\text{-}B).$$

The coherence of geometric $\Z$-algebra $\cA$ of period $n$ is
equivalent to some statement about t-structures. Namely, let
$(\tau_{\leq 0},\tau_{\geq 1})$ be a t-structure on
$D^b(\text{Mod-}\cA_{[1,n]})$ induced by the equivalence of
Theorem \ref{D(Geom)}. It can be shown that $\cA$ is coherent iff
the t-structure $(\tau_{\leq 0},\tau_{\geq 1})$ induces a
t-structure on $D^b(\text{mod-}_{finite}(\cA_{[1,n]})).$}
\end{remark}

\section{The $k$-points of noncommutative Grassmanians}
\label{k-points}

To discuss the $k$-points of noncommutative Grassmanians defined
above we should first relate the following two approaches to
noncommutative geometry.

The first one is to think of noncommutative stacks as of
$\Proj(\cA),$ where $\cA$ is a $\Z$-algebra. The special case of
graded algebras, i.e. $1$-periodic $\Z$-algebras is studied in
\cite{M}, \cite{V1}, \cite{V2}, \cite{AZ} and other papers. However,
it seems to be more reasonable to consider $\Z$-algebras. Note that
our noncommutative Grassmanians are naturally defined as $\Proj$ of
a $\Z$-algebra but not a graded algebra.

The other approach is to think of a noncommutative stacks as of
(equivalence classes of) presheaves of (small) groupoids $X$ on the
category $\Alg_k^{op}$ opposite to the category of unital
associative algebras. Morally the groupoid $X(A)$ should be thought
of as the groupoid of maps from the affine noncommutative scheme
$Sp(A)$ to $X.$ This approach is studied in \cite{Or} in the case of
sets (trivial groupoids). In this case we have the category of
quasi-coherent sheaves which is not always abelian (it always has
cokernels but may not admit kernels), and the structure sheaf.

In the second approach we obviously have the groupoid of $k$-points
$X(k).$

From this moment we assume that the $\Z$-algebra $\cA$ is positively
oriented but not necessarily connected. We will make an attempt to
define the presheaf of groupoids of morphisms $Sp(A)\to \Proj(\cA),$
$A\in \Alg_k.$ First note that a morphism $f:Sp(A)\to \Proj(\cA)$
must give a $k$-linear additive functor $f^*:\QMod(\cA)\to
\text{Mod-}A$ together with an isomorphism $f^*(\pi(P_0))\cong A.$
Moreover, $f^*$ must commute with colimits.

Notice that if $\cC$ is a $k$-linear abelian category with
infinite direct sums and $\cF:\cA\to \cC$ is a $k$-linear functor
then we have the tensor product functor

$$-\otimes_{\cA}\cF:\text{Mod-}\cA\to \cC$$ given by the formula

$$M\otimes_{\cA}\cF=\coker(b: M\otimes \cA\otimes\cF\to M\otimes \cF),$$
where $b=\mu_M\otimes \one_{\cA}-\one_{M}\otimes \mu_{\cF}$ (we
identify $\cF$ with $\bigoplus\limits_{i\in \Z} \cF(i)\in \cC$).
Clearly, the functor $-\otimes_{\cA}\cF$ commutes with colimits.
We denote by $\Tor_i^{\cA}(-,\cF)$ its left derived functors.

\begin{defi}Let $\cA$ be a positively oriented $\Z$-algebra
and let $\cC$ be a $k$-linear abelian category with infinite direct
sums together with a distinguished object $Y\in \cC.$

We denote by $\cG_1(\cA,\cC,Y)$ the groupoid of pairs
$(f^*,\theta),$ where $f^*:\QMod(\cA)\to \cC$ is a $k$-linear
functor commuting with colimits and $\theta:f^*(\pi(P_0))\to Y$ is
an isomorphism.

We denote by $\cG_2(\cA,\cC,Y)$ the groupoid of pairs
$(\cF,\sigma),$ where $\cF:\cA\to \cC$ is a $k$-linear functor such
that $\Tor_0^{\cA}(T,\cF)=\Tor_1^{\cA}(T,\cF)=0$ for each torsion
$\cA^{op}$-module $T,$ and $\sigma:\cF(0)\to Y$ is an
isomorphism.\end{defi}

\begin{theo}\label{G_1=G_2} Let $\cA$ be a positively oriented $\Z$-algebra
and let $\cC$ be a $k$-linear abelian category with infinite direct
sums together with a distinguished object $Y\in \cC.$ Then the
groupoids $\cG_1(\cA,\cC,Y)$ and $\cG_2(\cA,\cC,Y)$ are equivalent.
\end{theo}

\begin{proof} We define the functor $\Phi:\cG_1(\cA,\cC,Y)\to \cG_2(\cA,\cC,Y)$ as follows. Let
$(f^*,\theta)\in \cG_1(\cA,\cC,Y).$ The functor $\Phi(f^*):\cA\to
\cC $ is defined by the formulas
$$\Phi(f^*)(i)=f^*(\pi(P_i)),$$
and for $x\in \cA_{ij}$ $$\Phi(f^*)(x)=f^*(\pi(x)).$$ We claim that
the pair $(\Phi(f^*),\theta)$ is an object of $\cG_2(\cA,\cC,Y).$
Indeed, let $T\in \Tors(\cA).$ Since the sequence
$$\pi(T\otimes \cA\otimes\cA\otimes\cA)\to \pi(T\otimes \cA\otimes\cA)\to
\pi(T\otimes\cA)\to 0$$ is exact in $\QMod(\cA)$ it follows that
the sequence
$$T\otimes \cA\otimes\cA\otimes\Phi(f^*)(\cA)\to T\otimes \cA\otimes \Phi(f^*)(\cA)\to
T\otimes \Phi(f^*)(\cA)\to 0$$ is exact in $\cC,$ i.e.
$$\Tor_0^{\cA}(S_i,\Phi(f^*))=\Tor_0^{\cA}(S_i,\Phi(f^*))=0.$$
Thus, the functor $\Phi$ is defined on objects. It obviously
extends to morphisms.

Now we define the functor $\Psi:\cG_2(\cA,\cC,Y)\to
\cG_1(\cA,\cC,Y)$ as follows. Let $(\cF,\sigma)\in
\cG_2(\cA,\cC,Y).$ We claim that the formula
$$\Psi(\cF)(\pi(M))=M\otimes_{\cA} \cF$$
well defines a functor $\Psi(\cF):\QMod(\cA)\to \cC$ which is
right exact and commutes with infinite direct sums. Indeed it
follows from the condition
$$\Tor_0^{\cA}(T,\cF)=\Tor_1^{\cA}(T,\cF)=0$$ for torsion $\cA^{op}$-modules $T.$

Hence, the pair $(\Psi(\cF),\sigma)$ is an object of
$\cG_1(\cA,\cC,Y).$ This defines the functor $\Psi$ on objects and
it obviously extends to morphisms.

It is clear that the composition $\Phi\circ\Psi$ is isomorphic to
the identity functor. To see this for the composition
$\Psi\circ\Phi,$ it remains to note that each functor $f^*$ from the
pair in $\cG_1(\cA,\cC,Y)$ can be reconstructed from the functor
$\Phi(f^*)$ using exact sequences $$\pi(M\otimes \cA\otimes\cA)\to
\pi(M\otimes\cA)\to \pi(M)\to 0.$$
\end{proof}

Notice that it follows from the above theorem that each functor
$f^*:\QMod(\cA)\to \cC$ commuting with direct sums and right exact
has the right adjoint $f_*:\cC\to \QMod(\cA)$ given by the formula
$f_*(X)=\pi(\widetilde{f_*(X)})$
$$\widetilde{f_*(X)}(i)=\Hom_{\cC}(f^*(\pi(P_i)),X),$$
and for $\phi\in \Hom_{\cC}(f^*(\pi(P_j)),X),$ $x\in \cA_{ij}$
$$\phi\cdot x=\phi\cdot f^*(\pi(x)).$$
Indeed, this follows from the formula
$$f^*(\pi(M))=M\otimes_{\cA}\cF.$$

It is clear that $X_i(\cA)=(A\mapsto \cG_i(\cA,\text{Mod-}A,A)),$
$i=1,2,$ are presheaves of groupoids on the category $\Alg_k^{op},$
and the equivalence from the above theorem extends to the
equivalence of these presheaves.

However, not all functors $f^*$ commuting with colimits should come
from true morphisms $f:Sp(A)\to \Proj(\cA).$ Although a true
presheaf of groupoids should be defined as a full (small)
subpresheaf of $X_2(\cA).$ We are going to make an attempt in this
direction. Our motivation is the following Proposition.

\begin{prop}\label{Lf^*} Let $\cA$ be a $\Z$-algebra. Further, let $\cC$ be a $k$-linear abelian category with
infinite direct sums and with the distinguished object $Y.$ Let
$(f^*,\theta)\in \cG_1(\cA,\cC,Y)$ and $(\cF,\sigma)\in
\cG_2(\cA,\cC,Y)$ be objects which correspond to each other under
the equivalence of Theorem \ref{G_1=G_2}. The following conditions
are equivalent:

(i) there exists a left derived functor $\bL f^*:D^-(\QMod(\cA))\to
D^-(\cC),$ and $\bL^if^*(\pi(P_j))=0$ for $i\ne 0$ and all $j$;

(ii) we have $\Tor_i^{\cA}(T,\cF)=0$ for all $i\geq 0,$ $T\in
\Tors(\cA).$
\end{prop}

\begin{proof} Prove that (i) implies (ii). We have that the functor
$f^*\cong -\otimes \cF$ maps acyclic right bounded complexes of
direct sums of $\pi(P_i)$ to acyclic complexes. Applying this to the
projection of the free resolution of a torsion module $T,$ we obtain
that $\Tor_i^{\cA}(T,\cF)=0$ for $i\geq 0.$

Prove that (ii) implies (i). Since each object in $\QMod(\cA)$ can
be covered by a direct sum of $\pi(P_i),$ it suffices to prove that
$f^*$ maps right bounded acyclic complexes of direct sums of
$\pi(P_i)$ to acyclic complexes.

Since the kernel and the cokernel of the morphism $P_j\to
\omega\pi(P_j)$ are torsion, we have
$\Tor_i^{\cA}(\omega\pi(P_j),\cF)=0$ for $i>0.$ Further, if
$K^{\cdot}$ is a right bounded acyclic complex of direct sums of
$\pi(P_i)$ then $\omega(K^{\cdot})$ has torsion cohomology.
Therefore,
$$f^*(K^{\cdot})=\omega(K^{\cdot})\otimes_{\cA}\cF=
\omega(K^{\cdot})\stackrel{\bL}{\otimes}_{\cA}\cF,$$ and the last
complex is acyclic since $\Tor_i^{\cA}(T,\cF)=0$ for $i\geq 0,$
$T\in \Tors(\cA).$
\end{proof}

\begin{defi}  For each $\Z$-algebra $\cA$ we define the presheaf
$X_{\cA}$ of groupoids on the category $\Alg_{k}^{op}$ as follows.
It is a full subpresheaf of $X _2(\cA)$ and the groupoid
$X_{\cA}(A)\subset X_2(\cA)(A)$ consists of pairs $(\cF,\sigma)\in
X_2(\cA)(A)$ such that:

1) we have $\Tor_i^{\cA}(T,\cF)=0$ for all $i$ and $T\in
\Tors(\cA)$;

2) the $A^{op}$-modules $\cF(i)$ are flat.\end{defi}

It is clear that $X_{\cA}$ is indeed a subpresheaf of $X_2(\cA).$
For $f=(\cF,\sigma)\in X_{\cA}(A)$ we denote by $f^*:\QMod(\cA)\to
\text{Mod-}\cA$ the corresponding functor $-\otimes_{\cA} \cF.$ We
also regard the objects $f\in X_{\cA}(A)$ as maps from $Sp(A)$ to
$\Proj(\cA),$ where $Sp(A)$ is a noncommutative affine scheme
corresponding to $A.$

The following Lemma simplifies the complicated condition on
$\Tor_i.$ Recall the torsion $\cA^{op}$-modules $T_{p,q}$ from
Section \ref{Z-algebras}.

\begin{lemma}\label{sufffices_for_Tii} Let $\cA$ be a positively oriented connected
$\Z$-algebra and let $\cC$ be a $k$-linear abelian category with
infinite direct sums. Let $\cF:\cA\to \cC$ be a $k$-linear
functor.

Suppose that $\Tor_i^{\cA}(T_{j,j},\cF)=0$ for all $i$ and $j.$ Then
$\Tor_i^{\cA}(T,\cF)=0$ for all $i$ and all torsion
$\cA^{op}$-modules $T.$
\end{lemma}

\begin{proof} First note that if $T=T\one_j,$ then $T$ has a left resolution by direct sums
of $T_{j,j}.$ Hence, Lemma holds for such $T.$

Further, the torsion modules $T_{p,q},$ $p\leq q,$ have finite
filtrations with subquotients $T_m$ such that $T_m=T_m\one_m,$
$p\leq m\leq q.$ Hence $\Tor_i^{\cA}(T_{p,q},\cF)=0$ for all $i,$
and $p\leq q.$ Now Lemma follows from the observation that each
torsion module has a left resolution by direct sums of modules
$T_{pq}.$
\end{proof}

\begin{defi} We say that a positively oriented $\Z$-algebra $\cA$
satisfies the condition (**) if the following holds:

(i) the algebra $\cA$ is generated by its subspaces $\cA_0$ and
$\cA_1$;

(ii) for each $i\in \Z,$ the object $\pi(P_i)$ has a finite right
resolution by direct sums of $\pi(P_j)$ with $j>i.$
\end{defi}

The next Proposition motivates the condition (**).

\begin{prop}\label{**_is_good} Let $\cA$ be a $\Z$-algebra satisfying (**), and $A\in \Alg_k.$ Then for each
$f=(\cF,\sigma)\in X_{\cA}(A)$ we have $\Aut(f)=\{1\}.$
\end{prop}

\begin{proof} Let $g\in \Aut(f).$ Clearly, $g(0):\cF(0)\to \cF(0)$ is the identity morphism.
Further, for each $i\in \Z$ the surjection
$\cA_{i,i+1}\otimes\pi(P_i)\to \pi(P_{i+1})$ is mapped by $f^*$ to
the surjection $\cA_{i,i+1}\otimes \cF(i)\to \cF(i+1).$ Hence, we
obtain by increasing induction over $i$ that $g(i):\cF(i)\to \cF(i)$
is the identity for $i\geq 0.$

Finally, since $\cA$ satisfies (**), it follows from Proposition
\ref{Lf^*} that there exists an injection of the form $\pi(P_i)\to
\bigoplus\limits_{\alpha} \pi(P_{j_{\alpha}})$ with $j_{\alpha}>i$
which is mapped by $f^*$ to the injection $\cF(i)\to
\bigoplus\limits_{\alpha} \cF(j_{\alpha}).$ Hence, we obtain by
decreasing induction on $i$ that $g(i):\cF(i)\to \cF(i)$ is the
identity for all $i\in \Z.$
\end{proof}

Therefore, if $\cA$ satisfies (**) we may and will replace $X_{\cA}$
by the equivalent presheaf of trivial groupoids $\pi_0(X_{\cA}).$ It
is easily seen from the proof of the above Proposition that
$\pi_0(X_{\cA}(A))$ is a set. Thus, $X_{\cA}$ is a presheaf of sets.

Now we would like to compare our definition of morphisms from
$Sp(A)$ to $\Proj(\cA)$ with the morphisms from commutative
Noetherian $k$-schemes to commutative projective $k$-schemes.

Note that we can restrict the presheaf $X_{\cA}$ to the full
subcategory of $\Alg_k$ which consists of commutative Noetherian
$k$-algebras. Further, we can extend this restricted presheaf onto
the category $Noeth_k$ of all commutative Noetherian $k$-schemes.

\begin{defi} Let $\cA$ be a positively oriented $\Z$-algebra. We
define the presheaf $\X_{\cA}: Noeth_k^{op}\to {\bf Gpd}$ as
follows. The groupoid $\X_{\cA}(Y)$ is a full sub-groupoid of
$\cG_2(\cA,\QCoh(Y),\cO_Y)$ which consists of objects $(\cF,\sigma)$
such that the following conditions hold:

1) we have $\Tor_i^{\cA}(T,\cF)=0$ for all $i$ and $T\in
\Tors(\cA)$;

2) the sheaves $\cF(i)$ are locally flat.
\end{defi}

We also regard the objects of the groupoid $\X_{\cA}(Y)$ as maps
from $Y$ to $\Proj(\cA).$ The analogue of Proposition
\ref{**_is_good} obviously holds for Noetherian $k$-schemes instead
of associative algebras. For each commutative Noetherian $k$-scheme
$Y$ we denote by $Y^{\vee}: Noeth_k^{op}\to \Set$ the presheaf of
sets represented by $Y.$

Now let $Z\subset \PP(V)$ be a closed subscheme and let $\cA$ be a
$\Z$-algebra associated to its homogeneous coordinate ring
$\bigoplus\limits_{d\geq 0}S^d V^*/I.$

\begin{prop}\label{comm_situation} The $\Z$-algebra $\cA$ satisfies the condition (**). The presheafs of sets $Z^{\vee}$ and
$\X_{\cA}:Noeth_k^{op}\to \Set$ on the category $Noeth_k^{op}$ are
isomorphic.
\end{prop}

\begin{proof} Recall that the category $\QCoh(Z)$ is equivalent to
$\Proj(\cA)$ by Serre Theorem. The sheaves $\cO_Z(i)$ correspond
under this equivalence to $\pi(P_i).$

Let $f:Y\to Z$ be a morphism. Then the sheaves $f^*(\cO_Z(i))$ are
invertible and hence are locally flat. Further, $f^*$ maps acyclic
right bounded complexes of direct sums of $\cO_Z(i)$ to acyclic
complexes. Finally, we have an isomorphism $f^*(\cO_Z)\cong
\cO_{Y}.$ Thus, we have a morphism of presheaves $Z^{\vee}\to
\X_{\cA}.$

Conversely, let $Y\in Noeth_k$ and $g\in \X_{\cA}(Y).$ Notice that
for each $i\in \Z$ we have an acyclic Koszul complex on $\PP(V)$
twisted by $\cO_{\PP(V)}(i),$ and we can restrict it to $Z$:
$$0\to \cO_Z(i)=\Lambda^n V^*\otimes \cO_Z(i)\to \Lambda^{n-1} V^*\otimes\cO_Z(i+1)
\to\dots\to \cO_Z(i+n)\to 0.$$ In particular, the $\Z$-algebra
$\cA$ satisfies the condition (**).

Hence, we have acyclic complexes \begin{equation}\label{acyclic}
0\to g^*(\pi(P_i))=\Lambda^n V^*\otimes g^*(\pi(P_i))\to
\Lambda^{n-1}V^* \otimes g^*(\pi(P_{i+1})) \to\dots\to
g^*(\pi(P_{i+n}))\to 0.\end{equation} In particular, we have
surjections $V^*\otimes g^*(\pi(P_i))\to g^*(\pi(P_{i+1}))$ and
injections $g^*(\pi(P_i))\to V\otimes g^*(\pi(P_{i+1})).$ Since
$g^*(\pi(P_0))\cong \cO_Y,$ we obtain by increasing and decreasing
inductions on $i$ that all the sheaves $g^*(\pi(P_i))$ are coherent
and are non-zero on each connected component of $Y.$ Since they are
locally flat, they are locally free.

Further, put $\cL=g^*(\pi(P_1)).$ Clearly, $g$ can be reconstructed
from the surjective morphism $\phi:V^*\otimes A\to \cL$ using the
exact sequences $$\Lambda^2 V^*\otimes g^*(\pi(P_{i-1}))\to
V^*\otimes g^*(\pi(P_i))\to g^*(\pi(P_{i+1}))\to 0$$ and $$0\to
g^*(\pi(P_{i-1}))\to V\otimes g^*(\pi(P_i))\to \Lambda^2 V\otimes
g^*(\pi(P_{i+1}))$$ from commlexes (\ref{acyclic}). Suppose that
$\rank (\cL_{|Y_0})\geq 2$ on some connected component $Y_0\subset
Y.$ Then it is easy to see that the morphism
$\Lambda^{n-1}V^*\otimes \cO_{Y_0}\to \Lambda^{n-2}V^*\otimes
\cL_{|Y_0}$ is injective. Hence $g^*(\pi(P_{-1}))_{|Y_0}=0,$ a
contradiction. Thus, $\cL$ is an invertible sheaf.

The surjective morphism $\phi$ above defines a morphism
$\widetilde{g}:Y\to \PP(V).$ It follows that $\widetilde{g}^*\cong
g^*\iota^*,$ where $\iota:Z\to \PP(V)$ is the embedding. Hence,
$g^*(\pi(P_i))\cong \cL^{\otimes i}.$ Further, the induced morphism
of graded algebras
$$\bigoplus\limits_{d\geq 0}S^d V^*\to \bigoplus_{i\geq
0}H^0(\cL^{\otimes i})$$ passes through $\bigoplus\limits_{d\geq
0}S^d V^*/I,$ thus $\widetilde{g}$ passes through $Z,$ and we obtain
a morphism $Y\to Z.$ Hence we have a morphism of presheaves
$\X_{\cA}\to Z^{\vee}.$

The constructed morphisms of presheaves are inverse to each other.
Proposition is proved.
\end{proof}

Now we want to describe the $k$-points of noncommutative
Grassmanians.

\begin{lemma}\label{geom_**} Let $\cA$ be a geometric $\Z$-algebra of period $n.$ Then it
satisfies the condition (**).
\end{lemma}
\begin{proof} By definition, $\Alg_{\cA}$ is generated by $\cA_0$
and $\cA_1.$ Further, the projections $\pi(K_i)$ of Koszul complexes
are acyclic. The first non-zero term of $\pi(K_{i+n})$ equals to
$\pi(P_i)$ (since $\cA^!$ is Frobenious). Therefore, each $\pi(P_i)$
has the required right resolution.
\end{proof}

Denote by $\pr_r^{m,V}:D^{b}_{coh}(\PP(V))\to \T^{m,V}$ the functor
which is right adjoint to the inclusion $\iota:\T^{m,V}\to
D^{b}_{coh}(\PP(V)).$ The next theorem describes the $k$-points of
noncommutative Grassmanians and the objects in $\T^{m,V}$
corresponding to their structure sheaves under the equivalence of
Corollary \ref{T^m,V}.

\begin{theo}\label{all_k-points} The $k$-points of the noncommutative Grassmanian
$\NGr(m,V),$ i.e. the elements $f\in X_{\cA^{m,V}}(k),$ naturally
correspond to vector subspaces $W\subset V$ of dimension $1\leq \dim
W\leq m.$ Further, if $f$ corresponds to $W$ then $f_*(k)$ is a
perfect complex, and it corresponds to the object
$\pr_r^{m,V}(\cO_{\PP(W)})\in \T^{m,V}$ under the equivalence
$\Perf(\QMod(\cA))\cong \T^{m,V}.$
\end{theo}

\begin{proof} Let $f\in X_{\cA}(k).$ For each $i\in \Z$ we have the
natural acyclic complex
$$0\to f^*(\pi(P_i))\cong \cA^{m,V!^*}_{i+n,i}\otimes f^*(\pi(P_i))\to
\cA^{m,V!^*}_{i+n,i+1}\otimes f^*(\pi(P_{i+1}))\to \dots\to
f^*(\pi(P_{i+n}))\to 0.$$ In particular, we have surjective map
$V^*\otimes f^*(\pi(P_i))\to f^*(\pi(P_{i+1}))$ and injective map
$f^*(\pi(P_i))\to V\otimes f^*(\pi(P_{i+1})).$ Since
$f^*(\pi(P_0))\cong k,$ we obtain by increasing and decreasing
inductions over $i$ that all the spaces $f^*(\pi(P_i))$ are non-zero
and finite-dimensional. Further, using the exact sequences
$$\Lambda^2 V^*\otimes f^*(\pi(P_{i-1}))\to V^*\otimes
f^*(\pi(P_i))\to f^*(\pi(P_{i+1}))\to 0$$ and $$0\to
f^*(\pi(P_{i-1}))\to V\otimes f^*(\pi(P_i))\to \Lambda^2 V\otimes
f^*(\pi(P_{i+1})),$$ one can reconstruct $f$ from the injection
$f^*(\pi(P_{-1}))\hookrightarrow V\otimes f^*(\pi(P_0))\cong V.$
Thus, we can associate a non-zero vector subspace $W\subset V$ to
each $f\in X_{\cA}(k)$ and $f$ can be reconstructed from the
subspace $W.$ We will show that $W\subset V$ gives a $k$-point iff
$1\leq \dim W\leq m.$

First suppose that $\dim W>m$ and $W$ gives a $k$-point $f.$ Then
$$f^*(\pi(P_1))=\coker (\Lambda^{n-m-1}V\otimes W\to \Lambda^{n-m}V),$$
and the last space is zero since $\dim W>m.$ But $f^*(\pi(P_1))\ne
0,$ a contradiction.

Now let $1\leq \dim W=d\leq m.$ Let $S=\{E_i\}$ be a geometric helix
in $\T^{m,V}$ of period $n-m+1$ such that $E_j=\cO_{\PP(V)}(j)$ for
$m-n\leq j\leq 0.$ Then the endomorphism $\Z$-algebra of $S$ is
equivalent to $\cA^{m,V}.$ We define the functor $\cF:\cA\to
k\text{-Vect}$ by the formula
$$\cF(i)=\Hom(E_i,\cO_{\PP(W)})^{\vee}=\Hom(E_i,\pr_r^{m,V}(\cO_{\PP(W)}))^{\vee}.$$
We put $f=(\cF,id).$

Now we prove that $f\in X_{\cA^{m,V}}(k).$ By Lemma
\ref{sufffices_for_Tii}, it suffices to show that
$\Tor_i(S_j,\cF)=0$ for $i>0,$ $j\in \Z.$ Since the complexes

$$0\to E_i\cong \cA^{m,V!^*}_{i+n,i}\otimes
E_i\to \cA^{m,V!^*}_{i+n,i+1}\otimes E_{i+1}\to \dots\to
E_{i+n}\to 0$$ of objects in $D^b_{coh}(\PP(V))$ have zero
convolutions, it suffices to prove that

$$\Hom^i(E_j,\pr_r^{m,V}(\cO_{\PP(W)}))=\Hom^i(E_i,\cO_{\PP(W)})=0$$
for $i\ne 0,$ $j\in \Z.$ It is clear that this holds for $m-n\leq
j\leq 0.$ Further, we have a Koszul resolution of the sheaf
$\cO_{\PP(W)}$ on $\PP(V)$:

$$0\to\Lambda^{n-d}(V/W)^*\otimes \cO_{\PP(V)}(d-n)\to\dots\to(V/W)^*\otimes\cO_{\PP(V)}(-1)\to\cO_{\PP(V)}.$$
Thus, $\pr_r^{m,V}(\cO_{\PP(W)})$ is isomorphic to the complex
$$\dots\to0\to\Lambda^{n-m}(V/W)^*\otimes \cO_{\PP(V)}(m-n)\to\dots
\to(V/W)^*\otimes\cO_{\PP(V)}(-1)\to\cO_{\PP(V)}\to 0\to \dots.$$
Since the helix $\{E_i\}$ is geometric, we have
$\Hom^i(E_j,\pr_r^{m,V}(\cO_{\PP(W)}))=0$ for $i>0,$ $j\leq 0.$

Recall that $E_{j+n-m+1}=\Phi^{-1}(E_j),$ where $\Phi=F[m-n],$ and
$F$ is a Serre functor on $\T^{m,V}.$ We have that
\begin{equation}\label{Serre_inverse} F^{-1}(K^{\cdot})\cong\pr_l^{m,V}(K^{\cdot}\otimes
\cO_{\PP(V)}(n)[1-n]),\end{equation} where
$\pr_l^{m,V}:D^{b}_{coh}(\PP(V))\to \T^{m,V}$ is the functor which
is left adjoint to the inclusion $\iota:\T^{m,V}\to
D^{b}_{coh}(\PP(V)).$

\begin{lemma}\label{preserving} The functor $\pr_l^{m,V}:D^b_{coh}(\PP(V))\to \T^{m,V}$ maps $Ob(D^{\geq
i}_{coh}(\PP(V)))$ to $Ob(D^{\geq i-m+1}_{coh}(\PP(V))\cap
\T^{m,V}).$ The functor $\Phi^{-1}$ preserves $Ob(D^{\geq
i}_{coh}(\PP(V))\cap \T^{m,V}).$
\end{lemma}
\begin{proof} The second statement follows from the first one by the isomorphism
(\ref{Serre_inverse}). To prove the first statement, it suffices
to note that
$$\pr_l^{m,V}(X)\cong L_{\cO_{\PP(V)}(1)}\cdot\dots\cdot L_{\cO_{\PP(V)}(m-1)}(X)[m-1].$$
\end{proof}

Since $\pr_l^{m,V}(\cO(i))=0$ for $i=1,\dots,m-1,$ we have
$$\Phi^{-1}(\pr_r^{m,V}(\cO_{\PP(W)}))=\pr_l^{m,V}(\pr_r^{m,V}(\cO_{\PP(W)})\otimes \cO_{\PP(V)}(n))[1-m]=
\pr_l^{m,V}(\cO_{\PP(W)}(n))[1-m],$$ and the last object belongs to
$Ob(D^{\geq 0}_{coh}(\PP(V))\cap \T^{m,V})$ by Lemma
\ref{preserving}. Again by Lemma \ref{preserving} we have that
$\Phi^{-l}(\pr_r^{m,V}(\cO_{\PP(W)}))$ lies in $Ob(D^{\geq
0}_{coh}(\PP(V))\cap \T^{m,V})$ for $l>0.$ Thus, we have
$$\Hom^i(E_{j-(n-m+1)k},\pr_r^{m,V}(\cO_{\PP(W)}))=
\Hom^i(E_{j},\Phi^{-k}(\pr_r^{m,V}(\cO_{\PP(W)})))=0$$ for $i<0,$
$n-m\leq j\leq 0,$ and $k>0.$ Therefore,
$$\Hom^i(E_j,\pr_r^{m,V}(\cO_{\PP(W)}))=0$$ for $i\ne0,$ $j\leq 0.$

To prove the same for $j>0,$ note that
$\Hom^i(E_j,\pr_r^{m,V}(\cO_{\PP(W)}))$ is the $i$-th cohomology of
the complex
\begin{multline*}\dots\to0\to\Lambda^{n-m}(V/W)^*\otimes
\Hom^{n-m}(E_j,\cO_{\PP(V)}(m-n))\to\dots \\ \to(V/W)^*\otimes
\Hom^{n-m}(E_j,\cO_{\PP(V)}(-1))\to \Hom^{n-m}(E_j,
\cO_{\PP(V)})\to 0\to \dots,\end{multline*} where the left
non-zero term is in degree zero. This complex is dual to the
complex
\begin{multline*}\dots\to0\to\Hom(\cO_{\PP(V)},E_{j-n+m-1})\to (V/W)
\otimes \Hom(\cO_{\PP(V)}(-1),E_{j-n+m-1})\to\dots\\
\to\Lambda^{n-m}(V/W)\otimes \Hom(\cO_{\PP(V)}(m-n),E_{j-n+m-1})
\to 0\to \dots,\end{multline*} and the last one is isomorphic to
the complex
\begin{multline}\label{X_0} \dots\to0\to\Hom(E_{1-j},\cO_{\PP(V)}(m-n))\to (V/W)
\otimes \Hom(E_{1-j},\cO_{\PP(V)}(m-n+1))\to\dots\\
\to\Lambda^{n-m}(V/W)\otimes \Hom(E_{1-j},\cO_{\PP(V)}) \to 0\to
\dots.\end{multline} The complex (\ref{X_0}) calculates
$\Hom^i(E_{1-j},X_0),$ where $X_0$ is the complex
$$\dots0\to\cO_{\PP(V)}(m-n)\to (V/W)
\otimes \cO_{\PP(V)}(m-n+1)\to\dots \to\Lambda^{n-m}(V/W)\otimes
\cO_{\PP(V)} \to 0\to \dots.$$ Thus, it remains to show that
$\Hom^i(E_j, X_0)=0$ for $i<0,$ $j\leq 0.$

If we prove this for $n-m\leq j\leq 0,$ then the rest of the proof
will be analogous to the proof of the same vanishing for
$\pr_r^{m,V}(\cO_{\PP(W)})$ instead of $X_0.$ So let $i<0,$ $n-m\leq
j\leq 0.$ We have the chain of isomorphisms
\begin{multline*}\Hom^i(\cO_{\PP(V)}(j),X_0)\cong
\Hom^{n-m+i}(\pr_r^{m,V}(\cO_{\PP(W)}),\cO_{\PP(V)}(m-n-j))\cong\\
\Ext^{n-m+i}(\cO_{\PP(W)},\cO_{\PP(V)}(m-n-j))\cong\\
\Ext^{m-i-1}(\cO_{\PP(V)}(m-j),\cO_{\PP(W)})^{\vee}=H^{m-i-1}(\PP(W),\cO_{\PP(W)}(j-m))^{\vee},\end{multline*}
and the last space is zero since $m-i-1>d-1=\dim \PP(W).$

Thus, $f$ is indeed a $k$-point. Furthermore, we have that the
complex
$$\dots\to0\to\Lambda^{n-m}(V/W)^*\otimes \pi(P_{m-n})\to\dots
\to(V/W)^*\otimes\pi(P_{-1})\to\pi(P_0)\to 0\to \dots$$ is a
resolution of $f_*(k).$ Thus, $f_*(k)$ is a perfect complex and it
corresponds to the object $\pr_r^{m,V}(\cO_{\PP(W)})$ under the equivalence
$\Perf(\QMod(\cA^{m,V}))\cong \T^{m,V}.$ Theorem is proved.
\end{proof}

It turns out that the embedding of $k$-points
$\Gr(d,V)(k)\hookrightarrow \NGr(m,V)(k)$ for $1\leq d\leq m$ can be
extended to a morphism $\Gr(d,V)\to \NGr(m,V).$

\begin{prop}\label{Gr_to_NGr} Let $V$ be a finite-dimensional vector space and let
$1\leq d\leq m\leq \dim V=n.$ Then there exists a natural morphism
$f_{d,m,V}:\Gr(d,V)\to \NGr(m,V)$ such that the derived inverse
image functor $\bL f_{d,m,V}^*$ induces a full embedding
$$\Perf(\NGr(m,V))\to D^b_{coh}(\Gr(d,V)).$$
\end{prop}

\begin{proof} For each $W\in \Gr(d,V)$ denote by $f_{W}$ the corresponding $k$-point of
$\NGr(m,V).$ It is clear that there exist vector bundles $\cF(i)$ on
$\Gr(d,V)$ such that the fiber of $\cF(i)$ over the point
corresponding to $W$ is naturally identified with $f_W^*(\pi(P_i))$
(in particular, $\cF(-1)$ is a tautological bundle). So we have a
natural functor $\cF:\cA\to \QCoh(\Gr(m,V)).$ Also by Theorem
\ref{all_k-points} we have that the complexes of vector bundles
$$0\to \cF(i)\cong \cA^{m,V!^*}_{i+n,i}\otimes \cF(i)\to
\cA^{m,V!^*}_{i+n,i+1}\otimes \cF(i+1)\to \dots\to \cF(i+n)\to 0$$
are acyclic in the fibers over closed points (if the residue field
of a point is greater than $k$ we can make an extension of scalars).
Hence, these complexes are acyclic themselves. It follows from Lemma
\ref{sufffices_for_Tii} that the pair $(\cF,id)$ defines a map
$f_{d,m,V}:\Gr(d,V)\to \NGr(m,V).$

Further, for $m-n\leq j\leq 0$ we have that $\bL
f_{d,m,V}^*(\pi(P_j))=f_{d,m,V}^*(\pi(P_j))=S^{-j} E,$ where $E$ is
a tautological bundle. The collection
$(S^{n-m}E,\dots,E,\cO_{\Gr(d,V)})$ is a sub-collection of the full
strong exceptional collection on $\Gr(d,V)$ constructed by Kapranov
\cite{Ka}. Moreover, the functor $\bL f_{d,m,V}^*$ induces
isomorphisms
$$\Hom(\pi(P_i),\pi(P_j))\to \Hom(S^{-i}E,S^{-j}E)
$$ for $m-n\leq i\leq j\leq 0.$ Thus, the induced functor
$$\bL f_{d,m,V}^*:\Perf(\NGr(m,V))\to D^b_{coh}(\Gr(d,V))$$ is a full
embedding.
\end{proof}

Notice that the full embedding $\bL f_{1,m,V}^*:\Perf(\NGr(m,V))\to
D^b_{coh}(\PP(V))$ coincides with the composition of the equivalence
of Corollary \ref{T^m,V} with the tautological embedding
$\T^{m,V}\hookrightarrow D^b_{coh}(\PP(V)).$

\section{Completions of local rings of $k$-points}
\label{local_structure}

Let $X$ be a presheaf of sets on the category $\Alg_k^{op}$ of
noncommutative affine schemes. Let $x\in X(k)$ be a $k$-point.
Define the functor $F_{X,x}:\art\to \Set$ by the formula
$$F_{X,x}(\cR)=\{f\in X(\cR)| X(\iota)(f)=x\},$$
where $\iota:\cR\to k=\cR/m$ is the projection.

\begin{defi}\label{defi_local} Let $X$ be a presheaf of sets on the category $\Alg_k^{op}$ of
noncommutative affine schemes, and $x\in X(k)$ be a $k$-point. The
completion of the local ring $\widehat{\cO_x},$ if it exists, is
defined as a pro-artinian algebra which pro-represents the functor
$F_{X,x}:\art\to \Set.$
\end{defi}

We would like to describe the local rings of a $k$-point $x_W$ of
the noncommutative Grassmanian $\NGr(m,V)$ which correspond to a
subspace $W\subset V$ of dimension $m.$ Recall that by Lemma
\ref{geom_**} and by Proposition \ref{**_is_good} we have that
$X_{\cA^{m,V}}$ is (equivalent to) a presheaf of sets (trivial
groupoids). Thus, the above definition is applicable to
$X_{\cA^{m,V}}.$

\begin{prop}\label{REnd(O_W)} Let $W\subset V$ be a vector subspace, and $\dim V=n,$ $1\leq m=\dim W\leq n-1.$
Then the DG algebra $\bR\Hom^{\cdot}(\cO_{\PP(W)},\cO_{\PP(W)})$ is
formal, and the graded algebra
$\Ext^{\cdot}(\cO_{\PP(W)},\cO_{\PP(W)})$ is isomorphic to the
graded algebra
$$C^{W,V}=\bigoplus\limits_{d=0}^{n-m}\Lambda^d (V/W)\otimes S^d W^*.$$
\end{prop}
\begin{proof} Denote by $K_W^{\cdot}$ the Koszul resolution
$$0\to\Lambda^{n-m}(V/W)^*\otimes
\cO_{\PP(V)}(m-n)\to\dots\to(V/W)^*\otimes\cO_{\PP(V)}(-1)\to\cO_{\PP(V)}\to
0$$ of the sheaf $\cO_{\PP(W)}.$ Since $$\Ext^k
(\cO_{\PP(V)}(i),\cO_{\PP(V)}(j))=0$$ for $k>0,$ $m-n\leq i,j\leq
0,$ the DG algebra $\bR \Hom(\cO_{\PP(W)},\cO_{\PP(W)})$ is
quasi-isomorphic to the DG algebra
$\Hom^{\cdot}_{\cO_{\PP(V)}}(K_{W}^{\cdot},K_W^{\cdot}).$

Further, we also have that $$\Ext^k
(\cO_{\PP(V)}(i),\cO_{\PP(W)})=0$$ for $k>0,$ $m-n\leq i\leq 0.$
Thus, we have the chain of isomorphisms of graded vector spaces:

\begin{multline}\label{Ext=C^{W,V}}H^{\cdot}(\Hom^{\cdot}_{\cO_{\PP(V)}}(K_{W}^{\cdot},K_W^{\cdot}))\cong
H^{\cdot}(\Hom^{\cdot}_{\cO_{\PP(V)}}(K_{W}^{\cdot},\cO_{\PP(W)}))\cong\\
\cong \Hom^{\cdot}_{\cO_{\PP(V)}}(K_{W}^{\cdot},\cO_{\PP(W)})\cong
\bigoplus\limits_{d=0}^{n-m}\Lambda^d (V/W)\otimes S^d
W^*=C^{W,V}.\end{multline}

To prove the lemma, it suffices to construct a morphism of DG
algebras
$$
\varphi:C^{W,V}\to
\Hom^{\cdot}_{\cO_{\PP(V)}}(K_{W}^{\cdot},K_W^{\cdot}),
$$
which induces the identity in cohomology (under the isomorphisms
\ref{Ext=C^{W,V}}). To define $\varphi,$ one needs to define its
components
$$
\varphi_{d,i}:\Lambda^d(V/W)\otimes S^dW^*\to \Hom(\Lambda^{-i}(V/W)^*\otimes
\cO_{\PP(V)}(i), \Lambda^{-i-d}(V/W)^*\otimes \cO_{\PP(V)}(i+d))
$$
for $0\leq d\leq n-m,$ $m-n\leq i\leq -d.$ To do that, choose a
decomposition $V=W\oplus U.$ Then we have natural maps
$$\psi_{d,i}:\Lambda^dU\otimes S^dW^*\otimes
\Lambda^{-i}U^*\to \Lambda^{-i-d}U^*\otimes S^dV^*.$$ We define
$\varphi_{d,i}$ to be the maps corresponding to $\psi_{d,i}$ via
the isomorphisms
\begin{multline*}
\Hom(\Lambda^{-i}(V/W)^*\otimes \cO_{\PP(V)}(i),
\Lambda^{-i-d}(V/W)^*\otimes \cO_{\PP(V)}(i+d))\cong\\
\cong\Lambda^{-i}(V/W)\otimes \Lambda^{-i-d}(V/W)^*\otimes S^dV^*,
\end{multline*}
and
\begin{multline*}\Hom_k(\Lambda^d(V/W)\otimes
S^dW^*,\Lambda^{-i}(V/W)\otimes \Lambda^{-i-d}(V/W)^*\otimes
S^dV^*)\cong\\ \Hom_k(\Lambda^dU\otimes S^dW^*\otimes
\Lambda^{-i}U^*,\Lambda^{-i-d}U^*\otimes S^dV^*).
\end{multline*}

A straightforward checking shows that the map $\varphi$ with
components $\varphi_{d,i}$ satisfies the required properties.
\end{proof}

Note that the graded algebra
$C^{W,V}=\bigoplus\limits_{d=0}^{n-m}\Lambda^d (V/W)\otimes S^d
W^*$ is quadratic Koszul. Indeed, it coincides with "white" Manin
product $\bigoplus\limits_{d=0}^{n-m}\Lambda^d (V/W)\circ
\bigoplus\limits_{d\geq 0} S^d W^*,$ and according to \cite{PP}
the white product of quadratic Koszul algebras is again Koszul.

Thus, if we denote by $\hat{S}$ the dual of its (augmented) bar
construction, then we have that the projection $\hat{S}\to
H^0(\hat{S})$ is a quasi-isomorphism, and the algebra $H^0(\hat{S})$
is the completion of $C^{W,V!}$ with omitted grading. For
convenience we will write $C$ instead of $C^{W,V}.$ We denote by
$\mathfrak{A}$   the completion $\widehat{C^!}$ of the algebra $C^!$
with omitted grading.

\begin{theo}\label{local_rings} Let $W\subset V$ be a subspace of dimension $1\leq m\leq n-1.$
Let $x_{W}\in X_{\cA^{m,V}}(k)$ be the $k$-point of noncommutative
Grassmanian $\NGr(m,V)$ corresponding to the subspace $W\subset V.$
Then the algebra $\mathfrak{A}^{op}=\widehat{C^{!}}{}^{\; op}$
coincides with the completion of the local ring of the $k$-point
$x_{W}.$\end{theo}

\begin{proof} We will construct some morphism $u_W:Sp(
\mathfrak{A}^{op})\to \NGr(m,V)$ and then prove that it is the
universal one. For convenience, we will write $\cA$ instead of
$\cA^{m,V}.$

{\bf The construction of the morphism $u_W:Sp(\mathfrak{A}^{op})\to
\NGr(m,V).$} First we define an object $u_{W*}(\mathfrak{A})$ in the
category $\QMod(\cA),$ together with a morphism of algebras
$f:\mathfrak{A}^{op}\to \End(u_{W*}(\mathfrak{A})).$ Denote by
$M_W^{\cdot}$ the complex of $\cA^{op}-$modules
$$\dots\to 0\to \Lambda^{n-m}(V/W)^*\otimes P_{m-n}\to\dots\to (V/W)^*\otimes P_{-1}\to P_0\to
0\to \dots,$$ where $P_0$ is placed in degree zero. As in the proof
of Lemma \ref{REnd(O_W)}, the DG algebra
$\cB=\End_{\cA}(M_W^{\cdot})$ is quasi-isomorphic to $C;$ moreover,
each decomposition $V=W\oplus U$ gives a quasi-isomorphism $C\to
\cB.$ The DG algebra $\cB$ is naturally augmented: the augmentation
sends each $\phi\in \cB^0$ to its component
$\phi_0\in\End_{\cA}(P_0)=k.$ As usual, we have a natural element
$\alpha\in \cM\cC((B\bar{\cB})^*\otimes \cB).$ We put
$$M_W^{\cdot}\otimes_{\alpha} (B\bar{\cB})^*=(\cB\otimes_{\alpha}
(B\bar{\cB})^*)\otimes_{\cB} M_W^{\cdot}.$$ This is a DG
$(\cA\otimes (B\bar{\cB})^*)^{op}-$module. Thus,
$H^0(M_W^{\cdot}\otimes_{\alpha} (B\bar{\cB})^*)$ is an $(\cA\otimes
\mathfrak{A})^{op}-$module.

We put
$$
u_{W*}(\mathfrak{A})=\pi(H^0(M_W^{\cdot}\otimes_{\alpha}
(B\bar{\cB})^*))\in \QMod(\cA).
$$
The map $f:\mathfrak{A}^{op}\to
\End(u_{W*}(\mathfrak{A}))$ is the projection by $\pi$ of the
$\mathfrak{A}^{op}-$action on $H^0(M_W^{\cdot}\otimes_{\alpha}
(B\bar{\cB})^*).$

Define the functor $\cF:\cA\to \mathfrak{A}\text{-Mod}$ by the
formula
$$
\cF(i)=\Hom_{\mathfrak{A}^{op}}(\Hom_{\QMod(\cA)}(\pi(P_i),u_{W*}(\mathfrak{A})),\mathfrak{A}).
$$
Further, note that we have an isomorphism of DG
$(B\bar{\cB})^*-$modules
$$
(B\bar{\cB})^*\cong \Hom^{\cdot}_{\cA^{op}}(P_0,M_W^{\cdot}\otimes_{\alpha}(B\bar{\cB})^*).
$$

Passing to $H^0,$ we obtain an isomorphism
$$
\mathfrak{A}\to
\Hom_{\cA^{op}}(P_0,H^0(M_W^{\cdot}\otimes_{\alpha}(B\bar{\cB})^*)).
$$
Composing it with the projection by $\pi,$ we obtain the map
$$
\sigma':\mathfrak{A}\to
\Hom_{\QMod(\cA)}(\pi(P_0),u_{W*}(\mathfrak{A})).
$$

Applying the functor $\Hom_{\mathfrak{A}^{op}}(-,\mathfrak{A})$ to
the map $\sigma'$ we obtain the map
$$
\sigma:\cF(0)\to \mathfrak{A}.
$$

\begin{lemma}a) The map $\sigma'$ (and hence $\sigma$) is an isomorphism,
and $H^i(M_W^{\cdot}\otimes_{\alpha} (B\bar{\cB})^*)=0$ for $i\ne
0.$

b) The pair $(\cF,\sigma)$ defines an object of
$X_{\cA}(\mathfrak{A}^{op}),$ i.e. a morphism
$Sp(\mathfrak{A}^{op})\to \Proj(\cA)=\NGr(m,V).$
\end{lemma}

\begin{proof} a) Choose a decomposition $V=W\oplus U.$ As we
already mentioned above, such a decomposition gives a
quasi-isomorphism of DG algebras $C\to \cB.$ Composing the Koszul
dual morphism $(B\bar{\cB})^*\to \hat{S}$ with the projection
$\hat{S}\to H^0(\hat{S})=\mathfrak{A},$ we obtain the
quasi-isomorphism $\beta:(B\bar{\cB})^*\to \mathfrak{A}.$ Thus, we
may replace $M_W^{\cdot}\otimes_{\alpha}(B\bar{\cB})^*$ by
$M_W^{\cdot}\otimes_{\beta^*(\alpha)}\mathfrak{A}.$ The last object
is the complex of projective
$(\cA\otimes\mathfrak{A})^{op}-$modules. Furthermore, we have the
isomorphism of complexes of $\cA^{op}-$modules
$$(M_W^{\cdot}\otimes_{\beta^*(\alpha)}\mathfrak{A})\otimes_{\mathfrak{A}}k\cong M_W^{\cdot}.$$
Further, according to the proof of Theorem \ref{all_k-points}
$H^i(M_W^{\cdot})=0$ for $i\ne 0$ It follows that
$H^i(M_W^{\cdot}\otimes_{\beta^*(\alpha)}\mathfrak{A})=0$ for $i\ne
0.$

The space $\Hom_{\QMod(\cA)}(\pi(P_0),u_{W*}(\mathfrak{A}))$ is thus
the zeroth cohomology group of the complex
$\Hom^{\cdot}_{\QMod(\cA)}(\pi(P_0),\pi(M_W^{\cdot}\otimes_{\beta^*(\alpha)}\mathfrak{A})).$
It follows that $\sigma'$ is an isomorphism.

b) Firts prove that $\Ext^k(\pi(P_i),u_{W*}(\mathfrak{A}))=0$ for
$k>0,$ $i\in \Z,$ and $\Hom(\pi(P_i),u_{W*}(\mathfrak{A}))$ is a
free finitely generated $\mathfrak{A}^{op}-$module for $i\in \Z.$

For $i\leq 0,$ according to a), we have that
$\Ext^k(\pi(P_i),u_{W*}(\mathfrak{A}))$ is the $k-th$ cohomology of
the complex
$\Hom^{\cdot}(\pi(P_i),\pi(M_W^{\cdot}\otimes_{\beta^*(\alpha)}\mathfrak{A})).$
It is concentrated in non-positive degrees, hence
$\Ext^k(\pi(P_i),u_{W*}(\mathfrak{A}))=0$ for $k>0.$ Further, it is
bounded below complex of free finitely generated
$\mathfrak{A}^{op}-$modules with the only cohomology in degree zero.
Thus, this cohomology $\Hom(\pi(P_i),u_{W*}(\mathfrak{A}))$ is free
and finitely generated.

For $i>0,$ according to a), $\Ext^k(\pi(P_i),u_{W*}(\mathfrak{A}))$
is the $k-th$ cohomology of the complex
$\Ext^{n-m}(\pi(P_i),\pi(M_W^{\cdot}\otimes_{\beta^*(\alpha)}\mathfrak{A})[m-n]).$
This is a complex of free finitely generated
$\mathfrak{A}^{op}$-modules concentrated in degrees $0,1,\dots,n-m.$
If we multiply it by the left $\mathfrak{A}-$module $k,$ we will
obtain the complex $\Ext^{n-m}(\pi(P_i),\pi(M_W^{\cdot})[m-n]).$
This complex computes $\Ext^k(\pi(P_i),x_{W*}(k)).$ Thus, the only
cohomology of the source complex is in degree zero and is a free
finitely generated $\mathfrak{A}^{op}-$module.

Now we obtain that all $\cF(i)$ are free (and hence flat)
$\mathfrak{A}-$modules, and the complexes $$\dots\to 0\to
\cA_{k,k-n+m}^{!*}\otimes \cF(k-n+m)\to\dots\to
\cA_{k,k-1}^{!*}\otimes \cF(k-1)\to \cF(k)\to 0\to \dots$$ are
acyclic. By Lemma \ref{sufffices_for_Tii}, it follows that the pair
$(\cF,\sigma)$ defines an object of $X_{\cA}(\mathfrak{A}^{op})$
\end{proof}

We define $u_W$ as the pair $(\cF,\sigma).$ If $\pi:\mathfrak{A}\to
k$ is the projection, then by the very construction of $u_W,$ we
have $X_{\cA}(\pi)(u_W)=x_W.$

{\bf Universality.} Now we prove that the constructed morphism
$u_W$ is universal. More precisely, the pair $(\cF,\sigma)$ gives
the morphism of functors
$$\Phi:h_{\mathfrak{A}}\to F_{X_{\cA},x_W},$$
such that for each $f:\mathfrak{A}\to \cR,$
$\Phi(f)=X_{\cA}(f)(u_W).$ And we prove that $\Phi$ is an
isomorphism of functors.

In the proof of Theorem \ref{all_k-points} we have already seen that
each element $f\in X_{\cA}(k)$ is uniquely determined by the
injection $f^*(\pi(P_{-1}))\hookrightarrow V\otimes
f^*(\pi(P_0))=V.$ The same observation evidently holds for arbitrary
algebras $R$ instead of $k.$

Choose again a decomposition $V=W\oplus U.$ Choose some bases
$(e_1,\dots,e_m)$ of the vector space $W,$ and $(e_{m+1},\dots,e_n)$
of the vector space $U.$ Let $R$ be some complete local augmented
algebra with the (maximal) augmentation ideal $\m,$ and let $\pi:
R\to k$ be the projection. Let $f=(\cF'',\sigma'')\in X_{\cA}(R)$ be
some element such that $X_{\cA}(\pi)(f)=x_W.$ In particular, we have
a natural isomorphism $\cF''(-1)\otimes_A k\cong W.$ Since the
module $\cF''(-1)$ is flat, it is free. Let
$I:\cF''(-1)\hookrightarrow V\otimes R$ be the structure injection.
Clearly, there is a unique lift
$(\widetilde{e_1},\dots,\widetilde{e_m})$ onto $\cF''(-1)$ of the
basis $(e_1,\dots,e_m)$ such that
$$I(\widetilde{e_j})=e_j\otimes 1+\sum\limits_{i=m+1}^n e_i\otimes y_{ij},$$
where $y_{ij}\in \m.$ Thus, for each $(\cF'',\sigma'')$ as above
we have associated a matrix $(y_{ij}),$ $m+1\leq i\leq n,$ $1\leq
j\leq m,$ of elements in $\m.$ Moreover, $(\cF'',\sigma'')$ can be
reconstructed (up to a natural isomorphism) from this matrix.

If $g: R\to S$ is a morphism of complete local augmented algebras
and $(y_{ij})$ is the matrix associated to $f\in X_{\cA}(R),$ then
the matrix associated to $X_{\cA}(g)(f)$ equals to $(g(y_{ij})).$

In the case $R=\mathfrak{A}^{op},$ and $f=u_W=(\cF,\sigma),$ the
associated matrix is the following:
$$x_{ij}=e_i^*\otimes e_j\in (V/W)^*\otimes W\subset \mathfrak{A}^{op}.$$
The elements $x_{ij}$ are topological generators of the algebra
$\mathfrak{A}.$ The quadratic relations on them are the following:

\begin{equation}\label{rel1}[x_{ij},x_{lj}]=0\text{ for }m+1\leq i<l\leq n, 1\leq j\leq
m;\end{equation}
\begin{equation}\label{rel2}[x_{ij}+x_{ik},
x_{lj}+x_{lk}]=0\text{ for }m+1\leq i<l\leq n, 1\leq j<k\leq
m.\end{equation}

It follows that the morphism $\Phi$ is injective. To prove the
surjectivity and the theorem, it suffices to prove the following
Lemma.

\begin{lemma} Let $R$ be a local complete augmented algebra and
let $f=(\cF'',\sigma'')$ be as above. Let $(y_{ij}),$ $m+1\leq
i\leq n,$ $1\leq j\leq m,$ be the associated matrix, $y_{ij}\in
\m.$ Then the relations (\ref{rel1}), (\ref{rel2}) are satisfied
for $y_{ij}$ instead of $x_{ij}.$
\end{lemma}

\begin{proof}Using bases changes, we can
reduce the problem to the only relation $[y_{n-1,m},y_{n,m}]=0.$

We have $$\cF''(1)\cong \coker(\phi:\Lambda^{n-m-1}V\otimes
\cF''(-1)\to \Lambda^{n-m}V\otimes R).$$

Choose the basis of $(\widetilde{e_1},\dots,\widetilde{e_m})$ of the
$R^{op}$-module $\cF''(-1)$ as above. Then the explicit form of the
map $\phi$ is the following:

$$\phi((e_{i_1}\wedge\dots\wedge e_{i_{n-m-1}})\otimes \widetilde{e_j})=
(e_{i_1}\wedge\dots\wedge e_{i_{n-m-1}}\wedge e_j)\otimes
1+\sum\limits_{i=m+1}^{n}(e_{i_1}\wedge\dots\wedge
e_{i_{n-m-1}}\wedge e_i)\otimes y_{ij}.$$

It is clear that the vector space $\cF''(1)\otimes_R (R/\m)$ is
one-dimensional and generated by the projection of the element
$(e_{k+1}\wedge\dots\wedge e_n)\otimes 1.$ According to Nakayama
lemma and the condition on $\cF''(1)$ to be flat, the
$R^{op}-$module $\cF''(1)$ is freely generated by the projection of
the same element.

Finally, we notice that
\begin{multline*}
\phi((e_{k+1}\wedge\dots\wedge e_{n-1})\otimes
\widetilde{e_m}y_{n-1,m}+(e_{k+1}\wedge\dots\wedge e_{n-2}\wedge
e_n)\otimes \widetilde{e_m}y_{n,m}+(e_{k+1}\wedge\dots\wedge
e_{n-2}\wedge
e_m)\otimes\widetilde{e_m})\\=(e_{k+1}\wedge\dots\wedge
e_n)\otimes (y_{n,m}y_{n-1,m}-y_{n-1,m}y_{n,m}).
\end{multline*}
Therefore, the image of the RHS in $\cF''(1)$ is zero, and hence
$[y_{n-1,m},y_{n,m}]=0,$ q.e.d.
\end{proof}

Theorem is proved.
\end{proof}

\end{document}